\documentclass[a4paper, 12pt, twoside]{article}

\usepackage[a4paper, margin=1in]{geometry}
\usepackage{graphicx} 
\usepackage{amsmath,amsfonts,amsthm}
\usepackage[T1]{fontenc} 
\usepackage[english]{babel}

\usepackage{dsfont}
\usepackage{natbib}
\usepackage{color, colortbl}

\usepackage[hyperindex,breaklinks]{hyperref}

\newcommand{\Ben}{\begin{enumerate}}
\newcommand{\Een}{\end{enumerate}}
\newcommand{\Bit}{\begin{itemize}}
\newcommand{\Eit}{\end{itemize}}
\newcommand{\Beq}{\begin{equation}}
\newcommand{\Eeq}{\end{equation}}
\newcommand{\Ba}{\begin{align*}}
\newcommand{\Ea}{\end{align*}}
\newcommand{\Mb}{\mathbf}

\newtheorem{theorem}{Theorem}
\newtheorem{proposition}{Proposition}

\newtheorem{definition}{Definition}

\newtheorem{Rq}{Remark} 

\title{A central limit theorem for functions of stationary max-stable random fields on $\mathds{R}^d$}
\date{July 24, 2018}

\begin{document}

\author{Erwan Koch\footnote{EPFL, Chair of Statistics STAT, EPFL-SB-MATH-STAT, MA B1 433 (B\^atiment MA), Station 8, 1015 Lausanne, Switzerland. Email: erwan.koch@epfl.ch},~~ Cl\'ement Dombry\footnote{Universit\'e Bourgogne Franche-Comt\'e, Laboratoire de Math\'ematiques de Besan\c{c}on, UMR 6623, 16 route de Gray, 25030 Besan\c{c}on, France. Email: clement.dombry@univ-fcomte.fr},~~ Christian Y. Robert\footnote{ISFA Universit\'e Lyon 1, 50 Avenue Tony Garnier, 69366 Lyon cedex 07, France. \newline Email: christian.robert@univ-lyon1.fr}} 

\maketitle
 
\begin{abstract}
Max-stable random fields are very appropriate for the statistical modelling of spatial extremes. Hence, integrals of functions of max-stable random fields over a given region can play a key role in the assessment of the risk of natural disasters, meaning that it is relevant to improve our understanding of their probabilistic behaviour. For this purpose, in this paper, we propose a general central limit theorem for functions of stationary max-stable random fields on $\mathds{R}^d$. Then, we show that appropriate functions of the Brown-Resnick random field with a power variogram and of the Smith random field satisfy the central limit theorem. Another strong motivation for our work lies in the fact that central limit theorems for random fields on $\mathds{R}^d$ have been barely considered in the literature. As an application, we briefly show the usefulness of our results in a risk assessment context.

\medskip

\end{abstract}

\section{Introduction}

Max-stable random fields constitute an extension of extreme-value theory to the level of random fields \citep[in the case of stochastic processes, see, e.g.,][]{haan1984spectral, de2007extreme} and turn out to be fundamental for spatial extremes. Indeed, they are particularly well suited to model the temporal maxima of a given variable (for instance a meteorological variable) at all points in space since they arise as the pointwise maxima taken over an infinite number of appropriately rescaled independent and identically distributed (iid) random fields. Thus, appropriate functions of max-stable random fields can be adequate models for the costs triggered by extreme environmental events. Hence, normalized integrals on subsets of $\mathds{R}^2$ of functions of max-stable random fields and associated spatial risk measures \citep[see][]{koch2017spatial, koch2018SpatialRiskAxiomsGeneral} are useful for assessing the impact of natural disasters. The existence of a central limit theorem (CLT) for functions of max-stable random fields on $\mathds{R}^2$ would provide insights about the asymptotic probabilistic behaviour of the previously mentioned normalized integrals. Moreover, as explained in \cite{koch2017spatial, koch2018SpatialRiskAxiomsGeneral}, it is relevant to look at the evolution of spatial risk measures when the region over which they are applied becomes increasingly large; see the axiom of asymptotic spatial homogeneity of order $- \alpha$ in \cite{koch2017spatial, koch2018SpatialRiskAxiomsGeneral} which quantifies the rate of spatial diversification when the region becomes large. As shown in the latter paper, under relatively mild assumptions, asymptotic spatial homogeneity of order $-2$, $-1$ and $-1$ of spatial risk measures associated respectively with variance, Value-at-Risk and expected shortfall is satisfied when there is a CLT for the cost field. Finally, from a statistical viewpoint, the existence of a CLT allows to show the asymptotic normality of various estimators. For all these reasons, in this paper, we are interested in showing a CLT for functions of max-stable random fields on $\mathds{R}^d$. 

A CLT has already been proved in \cite{dombry2012strong} in the case of stationary\footnote{Throughout the paper, stationarity refers to strict stationarity.} max-infinitely divisible random fields on $\mathds{Z}^d$. Several CLTs for stochastic processes on $\mathds{Z}$ have been proposed under various (especially mixing) conditions  \citep[see, e.g.,][]{ibragimov1962some, gordin1969, ibragimov1975note}. Similarly, in the case of random fields on $\mathds{Z}^d$, many CLTs have been introduced under miscellaneous (especially mixing) conditions and in diverse contexts \citep[see, e.g.,][]{bolthausen1982central, chen1991uniform, nahapetian1992martingale, guyon1995random, perera1997geometry, dedecker1998central, el2013central}. For instance, the influential paper by \cite{bolthausen1982central} establishes a CLT for stationary random fields on $\mathds{Z}^d$ satisfying specific strong mixing conditions. However, the literature about CLTs for stochastic processes on $\mathds{R}$ or random fields on $\mathds{R}^d$ is, surprisingly, limited. This provides an additional strong motivation for our work. \cite{bulinski2010central} proposes a variant of the classical CLT where he considers a random field on $\mathds{R}^d$ but observed on a grid. His results involve two asymptotics at the same time: both the spatial domain and the grid resolution increase. The second type of asymptotics is known as infill asymptotics. To the best of our knowledge, only \cite{gorodetskii1984central, gorodetskii1987moment} has proposed a general CLT for strong mixing random fields on $\mathds{R}^d$. However, the strong mixing condition needed for this theorem to hold seems very difficult to check. Finally, CLTs for the indicator function of stationary random fields exceeding a given threshold have been obtained \citep[see, e.g.,][]{spodarev2014limit, bulinski2012central}. For more references about CLTs for random fields, we refer the reader for instance to \cite{ivanov1989statistical}. 

In this paper, we show a CLT for functions of stationary max-stable random fields on $\mathds{R}^d$. For the reason explained above, we could not use the results by \cite{gorodetskii1984central, gorodetskii1987moment} and, instead, we take advantage of the CLT by \cite{bolthausen1982central}. We basically extend Theorem 2.3 in \cite{dombry2012strong} from $\mathds{Z}^d$ to $\mathds{R}^d$ in the max-stable case, using the bound for the $\beta$-mixing coefficient established in that paper. Then, we show that appropriate functions of the Brown-Resnick random field with a power variogram and of the Smith random field satisfy the CLT. Finally, we briefly show the usefulness of our results in a risk assessment context.

The remainder of the paper is organized as follows. In Section \ref{Sec_Notations_CLT_Bolthausen}, we introduce some concepts about mixing as well as the previously mentioned CLT by \cite{bolthausen1982central} and give some insights about max-stable random fields. In Section \ref{Sec_CLT_Functions_Stationary_Maxstable_Processes}, we first establish our general CLT for functions of stationary max-stable random fields and then consider the cases of the Brown-Resnick and Smith random fields. We shortly illustrate our results in a risk assessment context in Section \ref{Sec_Application}. Finally, Section \ref{Sec_Conclusion} contains a brief summary as well as some perspectives. Throughout the paper, the elements belonging to $\mathds{R}^d$ or $\mathds{Z}^d$ for some $d \in \mathds{N} \backslash \{0 \}$ are denoted using bold symbols, whereas those in more general spaces are denoted using normal font. All proofs can be found in the Appendix.

\section{Some notations and concepts}
\label{Sec_Notations_CLT_Bolthausen}

In the following, ``$\bigvee$'' denotes the supremum when the latter is taken over a countable set. Moreover, $\overset{d}{=}$ and $\overset{d}{\to}$ stand for equality and convergence in distribution, respectively. Unless otherwise stated, in the case of random fields, distribution has to be understood as the set of all finite-dimensional distributions. Finally, let $\lambda$ denote the Lebesgue measure in $\mathds{R}^d$.

\subsection{Brief introduction to mixing and the central limit theorem by Bolthausen}

Let $(\Omega, \mathcal{F}, \mathds{P})$ be a probability space and $\mathcal{X}$ be a locally compact metric space. Let $\{ X(x) \}_{x \in \mathcal{X}}$ be a real-valued random field. For $S \subset \mathcal{X}$ a closed subset, we denote by $\mathcal{F}^X_S$ the $\sigma$-field generated by the random variables $\{ X(x): x \in S \}$. Moreover, for any $\mathcal{F}_1 \subset \mathcal{F}$, we denote by $\sigma(\mathcal{F}_1)$ the $\sigma$-field generated by $\mathcal{F}_1$. Let $S_1, S_2 \subset \mathcal{X}$ be disjoint closed subsets. The $\alpha$-mixing coefficient \citep[introduced by][]{rosenblatt1956central} between the $\sigma$-fields $\mathcal{F}^X_{S_1}$ and $\mathcal{F}^X_{S_2}$ is given by
$$ \alpha^X(S_1,S_2)=\sup\Big\{|\mathds{P}(A\cap B)-\mathds{P}(A)\mathds{P}(B)|: A\in \mathcal{F}^X_{S_1}, B\in \mathcal{F}^X_{S_2} \Big\}.$$
The $\beta$-mixing coefficient or absolute regularity coefficient \citep[attributed to Kolmogorov in][]{volkonskii1959some} between the $\sigma$-fields $\mathcal{F}^X_{S_1}$ and $\mathcal{F}^X_{S_2}$ is defined by
\Beq
\label{Eq_Def_2_Beta_Mixing}
\beta^X(S_1, S_2)= \frac{1}{2} \sup \left \{ \sum_{i=1}^{I} \sum_{j=1}^J | \mathds{P}(A_i \cap B_j) - \mathds{P}(A_i) \mathds{P}(B_j) | \right \},
\Eeq
where the supremum is taken over all partitions $\{ A_1,\dots , A_I \}$ and $\{ B_1, \dots , B_J \}$ of $\Omega$ with the
$A_i$'s in $\mathcal{F}^X_{S_1}$ and the $B_j$'s in $\mathcal{F}^X_{S_2}$. The inequality 
\Beq
\label{Eq_Majoration_Alphamixing_With_Betamixing}
\alpha^X(S_1, S_2) \leq \frac{1}{2} \beta^X(S_1, S_2), \quad \mbox{for all } S_1, S_2 \subset \mathcal{X},
\Eeq
will be very useful. 

We now present the previously mentioned CLT for stationary random fields in $\mathds{Z}^d$ due to \cite{bolthausen1982central} since it will play a key role in the proof of our general CLT for functions of stationary max-stable random fields (Theorem \ref{Th_General_Case}). Let $\mathcal{X}=\mathds{Z}^d$. For $\Mb{h}_1, \Mb{h}_2 \in \mathds{Z}^d$, let $d(\Mb{h}_1, \Mb{h}_2)=\max_{1 \leq i \leq d} \{ |\Mb{h}_1(i)-\Mb{h}_2(i) | \}$, where, for $\Mb{h} \in \mathds{Z}^d$, the $\Mb{h}(i)$, $1 \leq i \leq d$, are the components of $\Mb{h}$. If $S_1, S_2 \subset \mathds{Z}^d$, let $d(S_1, S_2)=\inf \{ d(\Mb{h}_1, \Mb{h}_2): \Mb{h}_1 \in S_1, \Mb{h}_2 \in S_2 \}$. For $\Lambda\subset\mathbb{Z}^d$, we note $|\Lambda|$ the number of elements in $\Lambda$ and $\partial\Lambda$ the set of elements $\Mb{h}_1\in\Lambda$ such that there exists $\Mb{h}_2 \notin \Lambda$ with $d(\Mb{h}_1,\Mb{h}_2)=1$.
For a real-valued stationary random field $\{ X(\Mb{h}) \}_{\Mb{h}\in\mathds{Z}^d}$, we note
\begin{equation}\label{eq:amix}
\alpha^X_{kl}(m)=\sup\Big\{\alpha^X(S_1,S_2):\ S_1, S_2 \subset \mathds{Z}^d,\ |S_1| \leq k, |S_2| \leq l,\ d(S_1,S_2)\geq m \Big\},
\end{equation}
defined for $m\geq 1$ and $k,l\in\mathds{N}\cup\{\infty\}$. Finally, let $\mathrm{Cov}$ denote the covariance.
\begin{theorem}\label{theo:BCLT}
Let $\{ X(\Mb{h}) \}_{\Mb{h}\in\mathds{Z}^d}$ be a real-valued stationary random field satisfying the following three conditions:
\begin{align}
&\sum_{m=1}^\infty m^{d-1}\alpha^X_{kl}(m)<\infty\quad  \mbox{for\ all\ \ } k\geq 1,\ l\geq 1\ \mbox{such\ that\ } k+l\leq 4;\label{Eq_TCL2}\\
&\alpha^X_{1\infty}(m)=o(m^{-d});\label{Eq_TCL1}\\
&\mathds{E}\big[|X(\Mb{0})|^{2+\delta}\big]<\infty \quad\mbox{and}\quad \sum_{m=1}^\infty m^{d-1} (\alpha^X_{11}(m) )^{\delta/(2+\delta)}<\infty\quad\mbox{for\ some}\ \delta>0.\label{Eq_TCL3}
\end{align}
\noindent Moreover, let $(\Lambda_n)_{n \in \mathds{N}}$ be a fixed sequence of finite subsets of $\mathds{Z}^d$ which increases to $\mathds{Z}^d$ and such that $\lim_{n\to\infty} |\partial \Lambda_n|/|\Lambda_n|=0$.
Then we have that $\sum_{\Mb{h}\in\mathds{Z}^d} | \mathrm{Cov}(X(\Mb{0}),X(\Mb{h})) | < \infty$ and, \\ if $\sigma^2=\sum_{\Mb{h}\in\mathds{Z}^d} \mathrm{Cov}(X(\Mb{0}),X(\Mb{h}))>0$, 
$$\frac{1}{|\Lambda_n|^{1/2}} \sum_{\Mb{h}\in\Lambda_n} (X(\Mb{h})-\mathds{E}[X(\Mb{h})]) \overset{d}{\rightarrow} \mathcal{N}(0,\sigma^2),\quad \mbox{as}\ n\to\infty,$$
where $\mathcal{N}(\mu, \sigma^2)$ denotes the normal distribution with expectation $\mu \in \mathds{R}$ and variance $\sigma^2>0$.
\end{theorem}

\subsection{Brief introduction to max-stable random fields}

\begin{definition}[Max-stable random field]
\label{Def_Maxstable_Processes}
A random field $\left \{ Z(\Mb{x}) \right \}_{\Mb{x} \in \mathds{R}^d}$ is said to be max-stable if there exist sequences of functions 
$\left( a_T(\Mb{x}), \Mb{x} \in \mathds{R}^d \right)_{T \geq 1}> 0$ and \\
$\left(  b_T(\Mb{x}), \Mb{x} \in \mathds{R}^d \right)_{T \geq 1} \in \mathds{R}$
such that, for all $T \geq 1$,
$$
\left \{ \frac{ \bigvee_{t=1}^T \left \{ Z_t(\Mb{x}) \right \}-b_T(\Mb{x} )}{a_T(\Mb{x} )}  \right \} _{\Mb{x} \in \mathds{R}^d}  \overset{d}{=} \{ Z(\Mb{x}) \}_{\Mb{x} \in \mathds{R}^d},
$$
where the $\{ Z_t(\Mb{x})\}_{\Mb{x} \in \mathds{R}^d }, t=1, \dots, T,$ are independent replications of $Z$.
\end{definition}
A max-stable random field is said to be simple if it has standard Fr\'echet margins, i.e., for all $\Mb{x} \in \mathds{R}^d$, $\mathds{P}(Z(\Mb{x}) < z)=\exp \left( -1/z \right), z>0$.

Now, let $\{ \tilde{T}_i(\Mb{x}) \}_{\Mb{x} \in \mathds{R}^d}, i=1, \dots, n,$ be independent replications of a random field $\{ \tilde{T}(\Mb{x}) \}_{\Mb{x} \in \mathds{R}^d}$. Let 
$\left( c_n(\Mb{x}), \Mb{x} \in \mathds{R}^d \right)_{n \geq 1} >0$ and $\left( d_n(\Mb{x}), \Mb{x} \in \mathds{R}^d \right)_{n \geq 1} \in \mathds{R}$ be sequences of functions. If there exists a non-degenerate random field $\{ G(\Mb{x}) \}_{\Mb{x} \in \mathds{R}^d}$ such that
$$
 \left \{ \frac{\bigvee_{i=1}^n \left \{ \tilde{T}_i(\Mb{x}) \right \} -d_n(\Mb{x})}{c_n(\Mb{x})} \right \}_{\Mb{x} \in \mathds{R}^d} \overset{d}{\rightarrow}  \left \{ G(\Mb{x}) \right \}_{\Mb{x} \in \mathds{R}^d}, \mbox{ for } n \to \infty,
$$
then $G$ is necessarily max-stable; see, e.g., \cite{haan1984spectral}. This explains the relevance of max-stable random fields in the modelling of spatial extremes. 

We know \citep[see, e.g.,][]{haan1984spectral} that any simple max-stable random field $Z$ can be written as
\Beq
\label{Eq_Spectral_Representation_Stochastic_Processes}
\left \{ Z(\Mb{x}) \right \}_{\Mb{x} \in \mathds{R}^d} \overset{d}{=} \left \{ \bigvee_{i=1}^{\infty} \{ U_i Y_i(\Mb{x}) \} \right \}_{\Mb{x} \in \mathds{R}^d},
\Eeq
where the $(U_i)_{i \geq 1}$ are the points of a Poisson point process on $(0, \infty)$ with intensity $u^{-2} \lambda(\mathrm{d}u)$ and the $Y_i, i\geq 1$, are independent replications of a random field $\{ Y(\Mb{x}) \}_{\Mb{x} \in \mathds{R}^d}$ such that, for all $\Mb{x} \in \mathds{R}^d$,
$\mathds{E}[Y(\Mb{x})]=1$. The random field $Y$ is not unique and is referred to as a spectral random field of $Z$. Conversely, any random field of the form \eqref{Eq_Spectral_Representation_Stochastic_Processes} is a simple max-stable random field. 

Below, we introduce the Brown-Resnick random field, defined in \cite{kabluchko2009stationary} as a generalization of the stochastic process introduced in \cite{brown1977extreme}. Recall that a random field $\{ W(\Mb{x}) \}_{\Mb{x} \in \mathds{R}^d}$ is said to have stationary increments if the distribution of the random field $\{ W(\Mb{x}+\Mb{x}_0)-W(\Mb{x}_0) \}_{\Mb{x} \in \mathds{R}^d}$ does not depend on the choice of $\Mb{x}_0 \in \mathds{R}^d$. Provided the increments of $W$ have a second moment, the variogram of $W$, $\gamma_W$, is defined by 
$$ \gamma_W(\Mb{x})=\mathrm{Var}(W(\Mb{x})-W(\Mb{0})), \quad \Mb{x} \in \mathds{R}^d,$$
where $\mathrm{Var}$ stands for the variance.
Moreover, for a random field $\{ W(\Mb{x}) \}_{\Mb{x} \in \mathds{R}^d}$ having a second moment, we introduce the function $\sigma_W$ defined by
$$ \sigma_W(\Mb{x})=\left[ \mathrm{Var}(W(\Mb{x})) \right]^{\frac{1}{2}}, \quad \Mb{x} \in \mathds{R}^d.$$ 
\begin{definition}[Brown-Resnick random field]
\label{Def_Brown_Resnick_Process}
Let $\{ W(\Mb{x}) \}_{\Mb{x} \in \mathds{R}^d}$ be a centred Gaussian random field with stationary increments and with variogram $\gamma_W$. Let us consider the random field $Y$ written as
$$ \{ Y(\Mb{x}) \}_{\Mb{x} \in \mathds{R}^d}=\left \{ \exp \left( W(\Mb{x})-\frac{\sigma^2_W(\Mb{x})}{2} \right) \right \}_{\Mb{x} \in \mathds{R}^d}.$$
Then the simple max-stable random field defined by \eqref{Eq_Spectral_Representation_Stochastic_Processes} with $Y$ is referred to as the Brown-Resnick random field associated with the variogram $\gamma_W$. In the following, we will call this field the Brown-Resnick random field built with $W$.
\end{definition}
It is worth noting that the Brown-Resnick random field is stationary \citep[see][Theorem 2]{kabluchko2009stationary} and that its distribution only depends on the variogram \citep[see][Proposition 11]{kabluchko2009stationary}.

We now present the Smith random field. Let $(U_i, \Mb{C}_i)_{i \geq 1}$ be the points of a Poisson point process on $(0,\infty) \times \mathds{R}^d$ with intensity measure $u^{-2} \lambda(\mathrm{d}u) \times \lambda(\mathrm{d}\Mb{c})$. Independently, let $f_i, i \geq 1$, be independent replicates of some non-negative random function $f$ on $\mathds{R}^d$ satisfying $\mathds{E} \left[ \int_{\mathds{R}^d} f(\Mb{x}) \ \lambda(\mathrm{d}\Mb{x}) \right]=1$. Then, the Mixed Moving Maxima (M3) random field
\Beq
\label{Eq_Mixed_Moving_Maxima_Representation}
\{ Z(\Mb{x}) \}_{\Mb{x} \in \mathds{R}^d}= \left \{ \bigvee_{i=1}^{\infty} \{ U_i f_i(\Mb{x}-\Mb{C}_i) \} \right \}_{\Mb{x} \in \mathds{R}^d}
\Eeq
is a stationary and simple max-stable random field. The so called Smith random field introduced by \cite{smith1990max} is a specific case of M3 random field.
\begin{definition}[Smith random field]
Let $Z$ be written as in \eqref{Eq_Mixed_Moving_Maxima_Representation} with $f$ being the density of a $d$-variate Gaussian random vector with mean $\Mb{0}$ and covariance matrix $\Sigma$. Then, the field $Z$ is referred to as the Smith random field with covariance matrix $\Sigma$.
\end{definition}

We conclude this section by giving some insights about a well-known dependence measure for max-stable random fields, the extremal coefficient. Let $\left \{ Z(\Mb{x}) \right \}_{\Mb{x} \in \mathds{R}^d}$ be a simple max-stable random field. For any compact $S \subset \mathds{R}^d$, the areal extremal coefficient of $Z$ for $S$, $\theta^Z(S)$, is defined by
\Beq
\label{Def_Extremal_Coefficient}
\mathds{P} \left( \sup_{\Mb{x} \in S} \{ Z(\Mb{x}) \} \leq z \right)=\exp \left( - \frac{\theta^Z(S)}{z} \right), \quad z>0.
\Eeq
It is easily shown that, for any compact $S \subset \mathds{R}^d$,
\Beq
\label{Eq_Property_Extremal_Coefficient}
\theta^Z(S)=\mathds{E} \left[ \sup_{\Mb{x} \in S} \{ Y(\Mb{x}) \} \right],
\Eeq
where $Y$ is a spectral random field of $Z$.

\section{A CLT for functions of stationary max-stable random fields on $\mathds{R}^d$}
\label{Sec_CLT_Functions_Stationary_Maxstable_Processes}

We start with some notations and definitions. Let $\|.\|$ stand for the Euclidean norm in $\mathds{R}^d$ and $'$ designate transposition. Moreover, for $\Mb{h}=(h_1, \dots, h_d)^{'} \in \mathds{Z}^d$, we adopt the notation $[\Mb{h}, \Mb{h}+1]=[h_1,h_1+1] \times \dots \times [h_d, h_d+1]$. We introduce $\mathds{S}= \left \{ \Mb{x} \in \mathds{R}^d: \| \Mb{x} \|=1 \right \}$, the unit sphere of $\mathds{R}^d$. Moreover, for two functions $f$ and $g$ from $\mathds{R}^d$ or $\mathds{Z}^d$ to $\mathds{R}$, the notations
$f(\Mb{h}) \underset{\| \Mb{h} \| \to \infty}{=} o(g(\Mb{h}))$ and $f(\Mb{h}) \underset{\| \Mb{h} \| \to \infty}{\sim} g(\Mb{h})$ mean that $\lim_{\| \Mb{h} \| \to \infty} f(\Mb{h})/g(\Mb{h})=0$ and $\lim_{\| \Mb{h} \| \to \infty} f(\Mb{h})/g(\Mb{h})=1$, respectively, where, for $L<\infty$, $\lim_{\| \Mb{h} \| \to \infty} f(\Mb{h})/g(\Mb{h})=L$ signifies $\lim_{h \to \infty} \sup_{\Mb{u} \in \mathds{S}} \left \{ | f(h \Mb{u})/g(h \Mb{u}) - L| \right \}=0$. Finally, $\lim_{\| \Mb{h} \| \to \infty} f(\Mb{h})=\infty$ must be understood as $\lim_{h \to \infty} \inf_{\Mb{u} \in \mathds{S}} \left \{ f(h \Mb{u}) \right \}=\infty$. For $V \subset \mathds{R}^d$ and $r>0$, we denote $N(V, r)=\{ \Mb{x} \in \mathds{R}^d: \mathrm{dist}(\Mb{x}, V) \leq r \}$, where $\mathrm{dist}$ designates the Euclidean distance. Furthermore, for $V \subset \mathds{R}^d$, we denote $\partial V$ the boundary of $V$. For a compact and convex set $V \subset \mathds{R}^d$, let $s(V)$ be the inradius of $V$, i.e., the largest $s>0$ such that $V$ contains a ball of radius $s$. Finally, let $\mathcal{B}(\mathds{R})$ and $\mathcal{B}((0, \infty))$ be the Borel $\sigma$-fields on $\mathds{R}$ and $(0, \infty)$, respectively.

We now present the concept of Van Hove sequence, which will play an important role in the following.
\begin{definition}
\label{Def_Van_Hove_Sequence}
A Van Hove sequence in $\mathds{R}^d$ is a sequence $( V_n )_{n \in \mathds{N}}$ of bounded measurable subsets of $\mathds{R}^d$ satisfying $V_n \uparrow \mathds{R}^d$, $\lim_{n \to \infty} \lambda(V_n)=\infty$, and $\lim_{n \to \infty} \lambda(N(\partial V_n, r) )/\lambda(V_n) =0, \mbox{ for all } r>0$. 
\end{definition}
Note that the assumption ``bounded'' does not always appear in the definition of a Van Hove sequence. It is worth mentioning that many sequences of bounded measurable subsets of $\mathds{R}^d$ are Van Hove sequences. For instance, any sequence $(V_n)_{n \in \mathds{N}}$ of compact convex subsets of $\mathds{R}^d$ such that $\lim_{n \to \infty} s(V_n)=\infty$ is a Van Hove sequence \citep[see, e.g.,][Lemma 3.11]{lenz2005ergodic}.

\subsection{CLT in the general case}

In the following, we say that a random field $\{ X(\Mb{x}) \}_{\Mb{x} \in \mathds{R}^d}$ such that, for all $\Mb{x} \in \mathds{R}^d$, $\mathds{E}[X(\Mb{x})^2]<\infty$, satisfies the CLT if
\Beq
\label{Eq_Convergence_Abs_Int_Cov_Thm}
\int_{\mathds{R}^d} | \mathrm{Cov}(X(\Mb{0}), X(\Mb{x})) | \  \lambda(\mathrm{d}\Mb{x}) < \infty,
\Eeq 
\Beq
\label{Eq_Def_Sigma_Square}
\sigma^2= \int_{\mathds{R}^d} \mathrm{Cov}(X(\Mb{0}), X(\Mb{x})) \  \lambda(\mathrm{d}\Mb{x})>0,
\Eeq
and, for any Van Hove sequence $(V_n)_{n \in \mathds{N}}$ in $\mathds{R}^d$,
\Beq
\label{Eq_Normality_Normalized_Integral}
\frac{1}{[\lambda(V_n)]^{1/2}} \int_{V_n} (X(\Mb{x})-\mathds{E}[X(\Mb{x})]) \  \lambda(\mathrm{d}\Mb{x}) \overset{d}{\rightarrow} \mathcal{N}(0, \sigma^2), \quad \mbox{as } n\to\infty.
\Eeq
The main result of this section, stated directly below, gives sufficient conditions such that a function of a stationary max-stable random field satisfies the CLT. 
\begin{theorem}
\label{Th_General_Case}
Let $\{ Z(\Mb{x}) \}_{\Mb{x} \in \mathds{R}^d}$ be a simple, stationary and sample-continuous max-stable random field and $F$ be a measurable function from $((0, \infty), \mathcal{B}((0,\infty)))$ to $(\mathds{R},\mathcal{B}(\mathds{R}))$ satisfying
\Beq
\label{Eq_Assumption_F}
\mathds{E}\left[ |F(Z(\Mb{0}))|^{2+\delta} \right]<\infty,
\Eeq
for some $\delta>0$ and
\Beq
\label{Eq_Condition_Positivity_Sigma2}
\sigma^2= \int_{\mathds{R}^d} \mathrm{Cov} \left(F(Z(\Mb{0})), F(Z(\Mb{x}))\right) \  \lambda(\mathrm{d}\Mb{x})>0.
\Eeq
Furthermore, assume that, for all $\Mb{h} \in \mathds{Z}^d$, 
\Beq
\label{Eq_Decay_Cubes_Y}
\mathds{E} \left[ \min \left \{ \sup_{\Mb{x} \in [0,1]^d} \{ Y(\Mb{x}) \},   \sup_{\Mb{x} \in [\Mb{h},\Mb{h}+1]}  \{ Y(\Mb{x}) \} \right \} \right] \leq K \| \Mb{h} \|^{-b},
\Eeq
for some $K>0$, $b > d \max  \left \{ 2, (2+\delta)/\delta \right \}$ and where $\{ Y(\Mb{x}) \}_{\Mb{x} \in \mathds{R}^d}$ is a spectral random field of $Z$ (see \eqref{Eq_Spectral_Representation_Stochastic_Processes}).
\noindent Then the random field $\{ X(\Mb{x}) \}_{\Mb{x} \in \mathds{R}^d}=\{ F(Z(\Mb{x})) \}_{\Mb{x} \in \mathds{R}^d}$ satisfies the CLT. 
\end{theorem}

It should be noted that this result constitutes, in the max-stable case, an extension of Theorem 2.3 in \cite{dombry2012strong} where the CLT for discrete random fields indexed by $\mathds{Z}^d$ is considered. Another connection with \cite{dombry2012strong} lies in the fact that we take advantage of the upper-bound for the $\beta$-mixing coefficient of simple and sample-continuous max-stable random fields established in that paper (Theorem 2.2).

We provide here the structure of the proof in order to convey some of the main ideas. For the detailed proof, we refer the reader to the Appendix. Without loss of generality, we assume that $\mathds{E}[X(\Mb{0})]=0$.
The proof is divided into three main parts. The first one is dedicated to the proof of \eqref{Eq_Convergence_Abs_Int_Cov_Thm}.
Then, the second and third ones show \eqref{Eq_Normality_Normalized_Integral}. Let $(V_n)_{n \in \mathds{N}}$ be a Van Hove sequence in $\mathds{R}^d$. In order to prove \eqref{Eq_Normality_Normalized_Integral}, we take advantage of the fact that, for all $n \in \mathds{N}$, 
\Beq
\label{Eq_Decomposition_Integral_Total}
\frac{1}{[\lambda(V_n)]^{\frac{1}{2}}} \int_{V_n} X(\Mb{x}) \ \lambda(\mathrm{d}\Mb{x})= I_{n,1} + I_{n,2},
\Eeq
where 
\Beq
\label{Eq_Def_In1_In2}
I_{n,1}=\frac{1}{[\lambda(V_n)]^{\frac{1}{2}}} \int_{A_n} X(\Mb{x}) \ \lambda(\mathrm{d}\Mb{x}) \quad \mbox{and} \quad I_{n,2} =\frac{1}{[\lambda(V_n)]^{\frac{1}{2}}} \int_{V_n \backslash A_n} X(\Mb{x}) \ \lambda(\mathrm{d}\Mb{x}),
\Eeq
with
$$A_n= \bigcup_{\Mb{h} \in \mathds{Z}^d:[\Mb{h}, \Mb{h}+1] \subset V_n} [\Mb{h}, \Mb{h}+1].$$
The second part of the proof is devoted to the study of $(I_{n,1})_{n \in \mathds{N}}$. For any $n \in \mathds{N}$, the domain of the related integral, $A_n$, consists of the union of all cubes $[\Mb{h}, \Mb{h}+1]$, for $\Mb{h} \in \mathds{Z}^d$, which are entirely included in $V_n$. As will be seen, considering such sets allows to deal with a random field on $\mathds{Z}^d$ instead of $\mathds{R}^d$. Thus, we show that the assumptions of Theorem \ref{theo:BCLT} (Bolthausen's theorem) are satisfied, obtaining that 
$$ I_{n,1} \overset{d}{\rightarrow} \mathcal{N}(0, \sigma^2), \mbox{ for } n \to \infty.$$
Finally, the third part concerns the study of $(I_{n,2})_{n \in \mathds{N}}$. For any $n \in \mathds{N}$, points belonging to the domain of the related integral, $V_n \backslash A_n$, are at a Euclidean distance not larger than $\sqrt{d}$ from the boundary of $V_n$. Hence, using the fact that $(V_n)_{n \in \mathds{N}}$ is a Van Hove sequence, we show that $\lim_{n \to \infty} \mathrm{Var}(I_{n,2})=0$, which allows to obtain \eqref{Eq_Normality_Normalized_Integral}. 

\begin{Rq}
It is worth mentioning that the left-hand side of \eqref{Eq_Decay_Cubes_Y} does not depend on the choice of the spectral random field $Y$. It only depends on the areal extremal coefficient function of $Z$. Indeed, the same computation as that leading to \eqref{Eq_Sum_Extremal_Coefficients} gives that 
\begin{align}
& \quad \ \mathds{E} \left[ \min \left \{ \sup_{\Mb{x} \in [0,1]^d} \{ Y(\Mb{x}) \},   \sup_{\Mb{x} \in [\Mb{h},\Mb{h}+1]}  \{ Y(\Mb{x}) \} \right \} \right] \nonumber
\\&=\theta^Z\left([ 0,1]^d \right)+\theta^Z\left([\Mb{h},\Mb{h}+1] \right)
-\theta^Z\left([ 0,1 ]^d \cup[\Mb{h},\Mb{h}+1] \right).
\label{Eq_Explanation_Main_Condition}
\end{align}
\end{Rq}

We now provide some insights about \eqref{Eq_Decay_Cubes_Y}, which is the main condition in Theorem \ref{Th_General_Case}. Using \eqref{Def_Extremal_Coefficient}, it follows from \eqref{Eq_Explanation_Main_Condition} that, for all $z >0$,
\begin{align}
& \quad \ \mathds{E} \left[ \min \left \{ \sup_{\Mb{x} \in [0,1]^d} \{ Y(\Mb{x}) \},   \sup_{\Mb{x} \in [\Mb{h},\Mb{h}+1]}  \{ Y(\Mb{x}) \} \right \} \right] \nonumber
\\& = - z \log \left( \mathds{P} \left( \sup_{\Mb{x} \in [0,1]^d} \{ Z(\Mb{x}) \} \leq z \right) \right) - z \log \left( \mathds{P} \left( \sup_{\Mb{x} \in [\Mb{h},\Mb{h}+1]} \{ Z(\Mb{x}) \} \leq z \right) \right) \nonumber
\\& \quad \ + z \log \left( \mathds{P} \left( \left \{  \sup_{\Mb{x} \in [0,1]^d} \{ Z(\Mb{x}) \} \leq z \right \} \bigcap \left \{ \sup_{\Mb{x} \in [\Mb{h},\Mb{h}+1]} \{ Z(\Mb{x}) \} \leq z \right \} \right) \right) \nonumber
\\& = z \log \left( \frac{\mathds{P} \left( \left \{  \sup_{\Mb{x} \in [0,1]^d} \{ Z(\Mb{x}) \} \leq z \right \} \bigcap \left \{ \sup_{\Mb{x} \in [\Mb{h},\Mb{h}+1]} \{ Z(\Mb{x}) \} \leq z \right \} \right)}{\mathds{P} \left( \sup_{\Mb{x} \in [0,1]^d} \{ Z(\Mb{x}) \} \leq z \right) \mathds{P} \left( \sup_{\Mb{x} \in [\Mb{h},\Mb{h}+1]} \{ Z(\Mb{x}) \} \leq z \right)} \right).
\label{Eq_Expression_Fraction_Probabilities}
\end{align}
Therefore, \eqref{Eq_Decay_Cubes_Y} implies that, for all $z >0$, 
\Beq
\label{Eq_Implication_1_Main_Condition}
\lim_{\| \Mb{h} \| \to \infty} \frac{\mathds{P} \left( \left \{  \sup_{\Mb{x} \in [0,1]^d} \{ Z(\Mb{x}) \} \leq z \right \} \bigcap \left \{ \sup_{\Mb{x} \in [\Mb{h},\Mb{h}+1]} \{ Z(\Mb{x}) \} \leq z \right \} \right)}{\mathds{P} \left( \sup_{\Mb{x} \in [0,1]^d} \{ Z(\Mb{x}) \} \leq z \right) \mathds{P} \left( \sup_{\Mb{x} \in [\Mb{h},\Mb{h}+1]} \{ Z(\Mb{x}) \} \leq z \right)}=1.
\Eeq
Consequently, \eqref{Eq_Decay_Cubes_Y} appears as a mixing condition. This is confirmed by the following fact. As can be seen in the proof of Theorem \ref{Th_General_Case}, \eqref{Eq_Decay_Cubes_Y} entails that, for all $\Mb{x} \in \mathds{R}^d$,
$$ 2-\theta^Z(\{ \Mb{0}, \Mb{x} \}) \leq K \| \Mb{x} \|^{-b}, $$
which gives that 
\Beq
\label{Eq_Condition_Mixing_Maxstable}
\lim_{\| \Mb{x} \| \to \infty} 2-\theta^Z(\{ \Mb{0}, \Mb{x} \}) =0.
\Eeq
From \cite{kabluchko2010ergodic}, Theorem 3.1, and the fact that this result can be extended to random fields on $\mathds{R}^d$, $d>1$ \citep[see, e.g.,][p.20]{DombryHDR2012}, we know that \eqref{Eq_Condition_Mixing_Maxstable} means that $Z$ is mixing.

Finally, we have, for all $z>0$, that 
\begin{align*}
& \quad \ \frac{\mathds{P} \left( \left \{  \sup_{\Mb{x} \in [0,1]^d} \{ Z(\Mb{x}) \} \leq z \right \} \bigcap \left \{ \sup_{\Mb{x} \in [\Mb{h},\Mb{h}+1]} \{ Z(\Mb{x}) \} \leq z \right \} \right)}{\mathds{P} \left( \sup_{\Mb{x} \in [0,1]^d} \{ Z(\Mb{x}) \} \leq z \right) \mathds{P} \left( \sup_{\Mb{x} \in [\Mb{h},\Mb{h}+1]} \{ Z(\Mb{x}) \} \leq z \right)} 
\\& = 1+ \frac{\mathds{E} \left[ \mathds{I}_{ \left \{  \sup_{\Mb{x} \in [0,1]^d} \{ Z(\Mb{x}) \} \leq z \right \} } \mathds{I}_{ \left \{ \sup_{\Mb{x} \in [\Mb{h},\Mb{h}+1]} \{ Z(\Mb{x}) \} \leq z \right \} } \right] - \mathds{E} \left[ \mathds{I}_{ \left \{  \sup_{\Mb{x} \in [0,1]^d} \{ Z(\Mb{x}) \} \leq z \right \} } \right] \mathds{E} \left [ \mathds{I}_{ \left \{ \sup_{\Mb{x} \in [\Mb{h},\Mb{h}+1]} \{ Z(\Mb{x}) \} \leq z \right \} } \right] }{\mathds{E} \left[ \mathds{I}_{ \left \{  \sup_{\Mb{x} \in [0,1]^d} \{ Z(\Mb{x}) \} \leq z \right \} } \right] \mathds{E} \left [ \mathds{I}_{ \left \{ \sup_{\Mb{x} \in [\Mb{h},\Mb{h}+1]} \{ Z(\Mb{x}) \} \leq z \right \} } \right]}.
\\& = 1+ \frac{ \mathrm{Cov} \left( \mathds{I}_{ \left \{  \sup_{\Mb{x} \in [0,1]^d} \{ Z(\Mb{x}) \} \leq z \right \} }, \mathds{I}_{ \left \{ \sup_{\Mb{x} \in [\Mb{h},\Mb{h}+1]} \{ Z(\Mb{x}) \} \leq z \right \} }  \right)}{\mathds{E} \left[ \mathds{I}_{ \left \{  \sup_{\Mb{x} \in [0,1]^d} \{ Z(\Mb{x}) \} \leq z \right \} } \right] \mathds{E} \left [ \mathds{I}_{ \left \{ \sup_{\Mb{x} \in [\Mb{h},\Mb{h}+1]} \{ Z(\Mb{x}) \} \leq z \right \} } \right]}.
\end{align*}
Therefore, it follows from \eqref{Eq_Expression_Fraction_Probabilities} and \eqref{Eq_Implication_1_Main_Condition} that
$$\mathds{E} \left[ \min \left \{ \sup_{\Mb{x} \in [0,1]^d} \{ Y(\Mb{x}) \},   \sup_{\Mb{x} \in [\Mb{h},\Mb{h}+1]}  \{ Y(\Mb{x}) \} \right \} \right] \underset{\| \Mb{h} \| \to \infty}{\sim} \frac{ \mathrm{Cov} \left( \mathds{I}_{ \left \{  \sup_{\Mb{x} \in [0,1]^d} \{ Z(\Mb{x}) \} \leq z \right \} }, \mathds{I}_{ \left \{ \sup_{\Mb{x} \in [\Mb{h},\Mb{h}+1]} \{ Z(\Mb{x}) \} \leq z \right \} }  \right)}{\mathds{E} \left[ \mathds{I}_{ \left \{  \sup_{\Mb{x} \in [0,1]^d} \{ Z(\Mb{x}) \} \leq z \right \} } \right] \mathds{E} \left [ \mathds{I}_{ \left \{ \sup_{\Mb{x} \in [\Mb{h},\Mb{h}+1]} \{ Z(\Mb{x}) \} \leq z \right \} } \right]}.$$
Hence, we deduce from \eqref{Eq_Decay_Cubes_Y} that, for all $z>0$,
$$ \lim_{\| \Mb{h} \| \to \infty} \mathrm{Cov} \left( \mathds{I}_{ \left \{  \sup_{\Mb{x} \in [0,1]^d} \{ Z(\Mb{x}) \} \leq z \right \} }, \mathds{I}_{ \left \{ \sup_{\Mb{x} \in [\Mb{h},\Mb{h}+1]} \{ Z(\Mb{x}) \} \leq z \right \} }  \right)=0.$$

In the next proposition, we provide conditions ensuring that \eqref{Eq_Condition_Positivity_Sigma2} is satisfied.
Before stating this result, we briefly recall the concept of association which plays an important role for max-stable random vectors.

\begin{definition}[Association]
\label{Def_Association}
A random vector $\Mb{R} \in \mathds{R}^q$, for $q \geq 1$, is said to be associated if $\mathrm{Cov}(g_1(\Mb{R}), g_2(\Mb{R})) \geq 0$ for all non-decreasing functions $g_i: \mathds{R}^q \to \mathds{R}$ such that $\mathds{E}[|g_i(\Mb{R})|] < \infty$ and $\mathds{E}[|g_1(\Mb{R}) g_2(\Mb{R})|] < \infty$ ($i=1,2$). Here, the term ``non-decreasing'' function must be understood in the following sense: for $\Mb{r}_1, \Mb{r}_2 \in \mathds{R}^q$, $\Mb{r}_1 \leq \Mb{r}_2$ implies $g_i(\Mb{r}_1) \leq g_i(\Mb{r}_2)$ ($i=1,2$), where $\Mb{r}_1 \leq \Mb{r}_2$ is a coordinatewise inequality.
\end{definition}

\begin{proposition}
\label{Prop_Condition_Sigma_Positive}
Let $\{ Z(\Mb{x}) \}_{\Mb{x} \in \mathds{R}^d}$ be a simple, stationary and sample-continuous max-stable random field.
For any function $F$ which is measurable from $((0, \infty), \mathcal{B}((0,\infty)))$ to $(\mathds{R},\mathcal{B}(\mathds{R}))$, satisfies \eqref{Eq_Assumption_F} and is moreover non-decreasing and non-constant, the random field $\{ X(\Mb{x}) \}_{\Mb{x} \in \mathds{R}^d}=\{ F(Z(\Mb{x})) \}_{\Mb{x} \in \mathds{R}^d}$ satisfies
$$ \sigma^2= \int_{\mathds{R}^d} \mathrm{Cov}(X(\Mb{0}), X(\Mb{x})) \  \lambda(\mathrm{d}\Mb{x})>0.$$
\end{proposition}

\subsection{CLT in the case of the Brown-Resnick and Smith random fields}

In this section, we show that if $\{ Z(\Mb{x}) \}_{\Mb{x} \in \mathds{R}^d}$ is the Smith or the Brown-Resnick random field with an appropriate variogram, then the random field $\{ F(Z(\Mb{x})) \}_{\Mb{x} \in \mathds{R}^d}$, where $F$ is as in Theorem \ref{Th_General_Case}, satisfies the CLT. In order to establish these results, we need the following proposition about the spectral random field of the Brown-Resnick random field. 

\begin{proposition}
\label{Prop_Upper_Bound_Expectation_Min_Brown_Resnick}
Let $\{ W(\Mb{x}) \}_{\Mb{x} \in \mathds{R}^d}$ be a centred Gaussian random field with stationary increments such that $W(\Mb{0})=0$. Moreover, assume that $W$ is a.s. bounded on $[0,1]^d$ and that, for $\Mb{h} \in \mathds{Z}^d$, 
\Beq
\label{Eq_Condition1_Variogram_TCL}
\sup_{\Mb{x} \in [0,1]^d} \{ \sigma_W^2(\Mb{h})-\sigma_W^2(\Mb{x}+\Mb{h}) \} \underset{\| \Mb{h} \| \to \infty}{=} o (\sigma_W^2(\Mb{h})),
\Eeq
and 
\Beq
\label{Eq_Condition2_Variogram_TCL}
\lim_{\| \Mb{h} \| \to \infty} \frac{\sigma_W^2(\Mb{h})}{\ln( \| \Mb{h} \|)} =\infty.
\Eeq
Then the random field $Y$ defined by
$$ \{ Y(\Mb{x}) \}_{\Mb{x} \in \mathds{R}^d}=\left \{ \exp \left( W(\Mb{x})-\frac{\sigma_W^2(\Mb{x})}{2} \right) \right \}_{\Mb{x} \in \mathds{R}^d}$$
satisfies Condition \eqref{Eq_Decay_Cubes_Y} for all $b>0$.
\end{proposition}

\begin{proposition}
\label{Prop_Variogram_Power}
Let $\{ W(\Mb{x}) \}_{\Mb{x} \in \mathds{R}^d}$ be a centred Gaussian random field with stationary increments such that $W(\Mb{0})=0$. Moreover, assume that  
\Beq
\label{Eq_Condition_Variogram_Power}
\sigma_W^2(\Mb{x})=\eta \| \Mb{x} \|^{\alpha}, \quad \Mb{x} \in \mathds{R}^d,
\Eeq
where $\eta>0$ and $\alpha \in (0, 2]$. Then $W$ is sample-continuous and the random field $Y$ defined by
$$ \{ Y(\Mb{x}) \}_{\Mb{x} \in \mathds{R}^d}=\left \{ \exp \left( W(\Mb{x})-\frac{\sigma_W^2(\Mb{x})}{2} \right) \right \}_{\Mb{x} \in \mathds{R}^d}$$
satisfies Condition \eqref{Eq_Decay_Cubes_Y} for all $b>0$.
\end{proposition}
\begin{Rq}
The field $W$ in Proposition \ref{Prop_Variogram_Power} is a (L\'evy) fractional Brownian field with Hurst parameter $\alpha/2 \in (0,1]$; see, e.g., \cite{Bierme2017random}, Section 1.2.3 and \cite{samorodnitsky1994stable}, Example 8.1.3. Its covariances are written
$$ \mathrm{Cov}(W(\Mb{x}), W(\Mb{y}))=\frac{\eta}{2} \left( \| \Mb{x} \|^{\alpha} + \| \Mb{y} \|^{\alpha} - \| \Mb{x}-\Mb{y} \|^{\alpha} \right), \quad \Mb{x}, \Mb{y} \in \mathds{R}^d.$$
\end{Rq}

Combining Theorem \ref{Th_General_Case} and Proposition \ref{Prop_Variogram_Power}, we directly obtain the following result.
\begin{theorem}
\label{Th_2_TCL_Brown_Resnick}
Let $\{ Z(\Mb{x}) \}_{\Mb{x} \in \mathds{R}^d}$ be the Brown-Resnick random field associated with the variogram $\gamma_W(\Mb{x})=\eta \| \Mb{x} \|^{\alpha}$, where $\eta >0$ and $\alpha \in (0,2]$, and $F$ be as in Theorem \ref{Th_General_Case}. Then $\{ F(Z(\Mb{x})) \}_{\Mb{x} \in \mathds{R}^d}$ satisfies the CLT.
\end{theorem}

\begin{Rq}
Using Proposition \ref{Prop_Upper_Bound_Expectation_Min_Brown_Resnick} and a very similar proof as that of Theorem \ref{Th_2_TCL_Brown_Resnick}, we obtain the following result. Let $\{ Z(\Mb{x}) \}_{\Mb{x} \in \mathds{R}^d}$ be the Brown-Resnick random field built with a random field $\{ W(\Mb{x}) \}_{\Mb{x} \in \mathds{R}^d}$ which is sample-continuous and whose variogram satisfies
$$\sup_{\Mb{x} \in [0,1]^d} \{ \gamma_W(\Mb{h})-\gamma_W(\Mb{x}+\Mb{h}) \}\underset{\| \Mb{h} \| \to \infty}{=}o(\gamma_W(\Mb{h})),$$
and 
$$ \lim_{\| \Mb{h} \| \to \infty} \frac{\gamma_W(\Mb{h})}{\ln( \| \Mb{h} \|)} =\infty.$$
Moreover, let $F$ be as in Theorem \ref{Th_General_Case}. Then $\{ F(Z(\Mb{x})) \}_{\Mb{x} \in \mathds{R}^d}$ satisfies the CLT.
\end{Rq}

Similarly as above, we obtain the following result for the Smith random field.

\begin{theorem}
\label{Th_TCL_Smith}
Let $\{ Z(\Mb{x}) \}_{\Mb{x} \in \mathds{R}^d}$ be the Smith random field with covariance matrix $\Sigma$ and $F$ be as in Theorem \ref{Th_General_Case}. Then the random field $\{ F(Z(\Mb{x})) \}_{\Mb{x} \in \mathds{R}^d}$ satisfies the CLT.
\end{theorem}

\section{Application to risk assessment}
\label{Sec_Application}

If the random field $\{ X(\Mb{x}) \}_{\Mb{x} \in \mathds{R}^2}$ is a model for the cost field (e.g., the economic cost or the cost for an insurance company) due to damaging events having a spatial extent (typically such as weather events), then, as detailed in \cite{koch2017spatial, koch2018SpatialRiskAxiomsGeneral}, the random variable 
$$ L_N(V_n)=\frac{1}{\lambda(V_n)} \int_{V_n} X(\Mb{x}) \ \lambda (\mathrm{d}\Mb{x}),$$
where $V_n \subset \mathds{R}^2$, is relevant for risk assessment. It can be interpreted as the loss per surface unit (or, less rigorously, as the loss per insurance policy) over the region $V_n$. If $X$ has a first moment and a constant expectation ($X$ is first-order stationary), then we have that
\begin{align}
\frac{1}{[\lambda(V_n)]^{1/2}} \int_{V_n} (X(\Mb{x})-\mathds{E}[X(\Mb{x})]) \  \lambda(\mathrm{d}\Mb{x}) &= \frac{1}{[\lambda(V_n)]^{1/2}} \int_{V_n} X(\Mb{x}) \  \lambda(\mathrm{d}\Mb{x}) - [\lambda(V_n)]^{1/2} \mathds{E}[X(\Mb{0})] \nonumber
\\& = [\lambda(V_n)]^{1/2} (L_N(V_n)-\mathds{E}[X(\Mb{0})]).
\label{Eq_Expression_Quantity_Function_L_N}
\end{align}
Hence, the asymptotic (when $V_n \uparrow \mathds{R}^2$) probabilistic behaviour of $L_N(V_n)$ can be derived from that of the left-hand side of \eqref{Eq_Expression_Quantity_Function_L_N}, quantity which appears in the CLT for the random field $X$, provided it exists. This explains the usefulness of a CLT in a risk assessment context.

As a short application, we now consider an insurance company called Ins. We assume that, during year $n$, Ins only covers the totality of the risk associated with a specified hazard over a whole continuous region, denoted by $V_n$ and referred to as the domain of Ins in the following. Moreover, let us assume that each insurance policy has a deductible $v>0$. Suppose that the process of the cost due to the mentioned hazard during a specified period (say one year) is given by a stationary and sample-continuous max-stable random field $\{ Z_G(\Mb{x}) \}_{\Mb{x} \in \mathds{R}^2}$ having standard Gumbel margins, i.e. such that, for all $\Mb{x} \in \mathds{R}^2$, $\mathds{P}(Z_G(\Mb{x}) \leq z)= \exp(-\exp(-z))$, $z \in \mathds{R}$. On the region $V_n$, this cost field is related to policies covered by Ins only. Thus, the normalized loss for Ins is given by
$$ L_{N}(V_n)= \frac{1}{\lambda(V_n)} \int_{V_n} (Z_G(\Mb{x})-v) \ \mathds{I}_{ \{Z_G(\Mb{x}) >v \} } \ \lambda(\mathrm{d} \Mb{x}).$$
Now, observe that $\{ Z_G(\Mb{x}) \}_{\Mb{x} \in \mathds{R}^2}=\{ \ln(Z(\Mb{x})) \}_{\Mb{x} \in \mathds{R}^2}$ for a simple, stationary and sample-continuous max-stable random field $\{ Z(\Mb{x}) \}_{\Mb{x} \in \mathds{R}^2}$. Hence, denoting $u = \exp(v)$, we have that 
$$ L_{N}(V_n)= \frac{1}{\lambda(V_n)} \int_{V_n} F(Z(\Mb{x})) \ \lambda(\mathrm{d} \Mb{x}),$$
where the function $F$ is defined by $F(z)=\ln \left( z/u \right) \ \mathds{I}_{ \{z>u\} }$, $z>0$. It is clear that $F$ is measurable from $((0, \infty), \mathcal{B}((0,\infty)))$ to $(\mathds{R},\mathcal{B}(\mathds{R}))$. Moreover, by construction, the random field $\{ \ln ( Z(\Mb{x}) ) \}_{\Mb{x} \in \mathds{R}^2}$ has Gumbel margins. Hence, for all $\delta>0$, we have that
$\mathds{E}\left[ |F(Z(\Mb{0}))|^{2+\delta} \right]< \infty$. In addition, $F$ is non-constant and non-decreasing and, thus, we deduce from Proposition \ref{Prop_Condition_Sigma_Positive} that \eqref{Eq_Condition_Positivity_Sigma2} is satisfied. Let us choose a $\delta>0$ and assume that $Z$ satisfies \eqref{Eq_Decay_Cubes_Y}. Furthermore, assume that the sequence of domains (over the years) of Ins, $(V_n)_{n \in \mathds{N}}$, is a Van Hove sequence in $\mathds{R}^2$. Then, applying Theorem \ref{Th_General_Case} and using \eqref{Eq_Expression_Quantity_Function_L_N}, we obtain that 
\Beq
\label{TCL_Application}
[\lambda(V_n)]^{1/2} (L_{N}(V_n)-\mathds{E}[F(Z(\Mb{0}))]) \overset{d}{\rightarrow} \mathcal{N}(0, \sigma^2), \quad \mbox{as } n\to\infty.
\Eeq
This gives the asymptotic probabilistic behaviour of the normalized loss suffered by Ins. If the sequence $(V_n)_{n \in \mathds{N}}$ is constant (i.e. if Ins does not plan to extend its domain) but the region $V_n$ is large enough, \eqref{TCL_Application} provides an approximation of the distribution of $L_{N}(V_n)$:
$$L_N(V_n) \approx \mathcal{N} \left(\mathds{E}[F(Z(\Mb{0}))], \frac{\sigma^2}{\lambda(V_n)} \right),$$
where $\approx$ means ``approximately follows''.
If $V_n$ increases in the Van Hove sense (e.g., if Ins extends its domain), \eqref{TCL_Application} for instance gives insights about how the Value-at-Risk of $L_N(V_n)$ evolves. This is related to the axiom of asymptotic spatial homogeneity of order $- \alpha$, see \cite{koch2017spatial, koch2018SpatialRiskAxiomsGeneral}. 

\section{Conclusion}
\label{Sec_Conclusion}

We have shown a general CLT for functions of stationary max-stable random fields on $\mathds{R}^d$. Moreover, we have seen that appropriate functions of the Brown-Resnick random field with a power variogram and the Smith random field satisfy the CLT. As briefly discussed, such results can be useful in a risk assessment context. Moreover, this paper proposes a new contribution to the limited literature about CLT for random fields on $\mathds{R}^d$. Future work might consist in relaxing the sample-continuity and stationarity assumptions on the max-stable random field $Z$ as well as letting the function $F$ depend on the location $\Mb{x}$ (with the notations of Theorem 2). 

\section*{Acknowledgements}
Erwan Koch would like to thank RiskLab at ETH Zurich and the Swiss Finance Institute (the project leading to this paper started as he was at ETH Zurich) as well as the Swiss National Science Foundation (project number 200021\_178824) for financial support. Furthermore, financial support by the Bourgogne-Franche-Comt\'e region (C. Dombry, grant OPE-2017-0068) is gratefully acknowledged. Finally, the authors would like to thank the Associate Editor and the referees for insightful suggestions.

\appendix
\section{Appendix: Proofs}

\subsection{For Theorem \ref{Th_General_Case}}

\begin{proof}
\noindent \textbf{Part 1: Proof of \eqref{Eq_Convergence_Abs_Int_Cov_Thm}}

\medskip
\medskip

Using \eqref{Eq_Assumption_F} and the stationarity of $X$ (by stationarity of $Z$), we have, for all $\Mb{x} \in \mathds{R}^d$, $\mathds{E} \left[ |X(\Mb{x})|^{2+\delta} \right]<\infty$. Hence, Davydov's inequality (\citeauthor{davydov1968convergence}, \citeyear{davydov1968convergence}, Equation (2.2); \citeauthor{ivanov1989statistical}, \citeyear{ivanov1989statistical}, Lemma 1.6.2) gives that
\Beq
\label{Eq_Davydov_Inequality}
| \mathrm{Cov}(X(\Mb{0}), X(\Mb{x})) | \leq 10 \left[ \alpha^X(\{ \Mb{0} \}, \{ \Mb{x} \}) \right]^{\frac{\delta}{2+\delta}} \left( \mathds{E} \left[ |X(\Mb{0})|^{2+\delta} \right] \right)^{\frac{1}{2+\delta}} \left( \mathds{E} \left[ |X(\Mb{x})|^{2+\delta} \right] \right)^{\frac{1}{2+\delta}}.
\Eeq
For all $\Mb{x} \in \mathds{R}^d$, since $F$ is Borel-measurable, $X(\Mb{x})=F(Z(\Mb{x}))$ is $\mathcal{F}^Z_{ \{ \Mb{x} \} }$-measurable. Hence, $\mathcal{F}^X_{ \{ \Mb{x} \} } \subset \mathcal{F}^Z_{ \{ \Mb{x} \} }$, which gives that, for all $\Mb{x} \in \mathds{R}^d$,
\Beq
\label{Eq_Majoration_Alpha_Mixing_X_With_Alpha_Mixing_Z}
\alpha^X \left( \{ \Mb{0} \}, \{ \Mb{x} \} \right) \leq \alpha^Z \left( \{ \Mb{0} \}, \{ \Mb{x} \} \right).
\Eeq
Moreover, using \eqref{Eq_Majoration_Alphamixing_With_Betamixing} and Corollary 2.2 in \cite{dombry2012strong}, we obtain that, for all $\Mb{x} \in \mathds{R}^d$,
\Beq
\label{Eq_Majoration_Alpha_Mixing_Z}
\alpha^Z \left( \{ \Mb{0} \}, \{ \Mb{x} \} \right) \leq 2 [2-\theta^Z(\{ \Mb{0}, \Mb{x} \})].
\Eeq
Therefore, the combination of \eqref{Eq_Majoration_Alpha_Mixing_X_With_Alpha_Mixing_Z} and \eqref{Eq_Majoration_Alpha_Mixing_Z} gives that
$$ \alpha^X \left( \{ \Mb{0} \}, \{ \Mb{x} \} \right) \leq 2 [2-\theta^Z(\{ \Mb{0}, \Mb{x} \})].$$
Hence, \eqref{Eq_Davydov_Inequality} and the stationarity of $X$ give that
\Beq
\label{Eq_Reformulation_Davydov_Inequality}
| \mathrm{Cov}(X(\Mb{0}), X(\Mb{x})) | \leq 2^{\frac{\delta}{2+\delta}} 10 \left( \mathds{E} \left[ |X(\Mb{0})|^{2+\delta} \right] \right)^{\frac{2}{2+\delta}} [2-\theta^Z(\{ \Mb{0}, \Mb{x} \})]^{\frac{\delta}{2+\delta}}.
\Eeq
Now, it follows from \eqref{Eq_Property_Extremal_Coefficient} that
$\theta^Z(\{ \Mb{0}, \Mb{x} \})= \mathds{E}[\max \{ Y(\Mb{0}), Y(\Mb{x}) \}]$. Thus, using the facts that, for all $\Mb{x} \in \mathds{R}^2$, $\mathds{E}[Y(\Mb{x})]=1$, and, for all $a, b \in \mathds{R}, a+b-\max \{a,b\}=\min\{a,b\}$, as well as the linearity of the expectation and \eqref{Eq_Decay_Cubes_Y}, we have that
\begin{align}
2-\theta^Z(\{ \Mb{0}, \Mb{x} \}) &=\mathds{E}[Y(\Mb{0})+Y(\Mb{x})-\max \{ Y(\Mb{0}), Y(\Mb{x}) \}] \nonumber
\\& =\mathds{E}[ \min \{ Y(\Mb{0}), Y(\Mb{x}) \}] \nonumber
\\& \leq \mathds{E} \left[ \min \left \{ \sup_{\Mb{y} \in [0,1]^d} \{ Y(\Mb{y}) \},   \sup_{\Mb{y} \in [\Mb{x},\Mb{x}+1]}  \{ Y(\Mb{y}) \} \right \} \right] \nonumber
\\& \leq K \| \Mb{x} \|^{-b}. \nonumber
\end{align}
As for all $\Mb{x} \in \mathds{R}^2$, $\theta^Z(\{ \Mb{0}, \Mb{x} \}) \leq 2$, and $b>(2+\delta)/\delta$, this directly implies that 
$$ \int_{\mathds{R}^d} \left[ 2-\theta^Z(\{ \Mb{0}, \Mb{x} \}) \right]^{\frac{\delta}{2+\delta}} \  \lambda(\mathrm{d}\Mb{x}) < \infty.$$ Finally, we obtain, using \eqref{Eq_Assumption_F} and \eqref{Eq_Reformulation_Davydov_Inequality}, that
\Beq
\label{Eq_Convergence_Int_Abs_Cov}
\int_{\mathds{R}^d} | \mathrm{Cov}(X(\Mb{0}), X(\Mb{x})) | \  \lambda(\mathrm{d}\Mb{x}) = K_1 < \infty,
\Eeq
which shows \eqref{Eq_Convergence_Abs_Int_Cov_Thm}. 

\medskip
\medskip

\noindent \textbf{Part 2: Study of $(I_{n,1})_{n \in \mathds{N}}$}

\medskip
\medskip

Introducing the random field
\Beq
\label{Eq_Def_Xtilde}
\tilde X(\Mb{h})=\int_{[\Mb{h}, \Mb{h}+1]} X(\Mb{x}) \ \lambda(\mathrm{d}\Mb{x}),\quad \Mb{h} \in \mathds{Z}^d,
\Eeq
we have, for all $n \in \mathds{N}$, that
\Beq
\label{Eq_Integral_Sum_Representation}
\int_{A_n} X(\Mb{x}) \  \lambda(\mathrm{d}\Mb{x})=\sum_{\Mb{h} \in \Lambda_n} \tilde X(\Mb{h}),
\Eeq
where $\Lambda_n=\{ \Mb{h} \in \mathds{Z}^d: [\Mb{h}, \Mb{h}+1] \subset A_n \}$. We will now show that the random field $\left \{ \tilde{X}(\Mb{h}) \right \}_{\Mb{h} \in \mathds{Z}^d}$ and the sequence $(\Lambda_n)_{n \in \mathds{N}}$ satisfy the assumptions of Theorem \ref{theo:BCLT} (Bolthausen's theorem). First, note that $\tilde{X}$ is stationary since $X$ is stationary.

\medskip

As already mentioned, since $F$ is Borel-measurable, we have, for all $\Mb{x} \in \mathds{R}^d$, that $X(\Mb{x})=F(Z(\Mb{x}))$ is $\mathcal{F}^Z_{ \{ \Mb{x} \} }$-measurable. Moreover, we know that the integral of measurable functions is again measurable. Hence, we have that $\tilde{X}(\Mb{h})$ is $\mathcal{F}_{[\Mb{h}, \Mb{h}+1]}^{Z}$-measurable. This directly gives that, for all $\Mb{h} \in \mathds{Z}^d$, 
\Beq
\label{Eq_Inclusion_SigmafieldXk_SigmafieldZ}
\mathcal{F}^{\tilde{X}}_{ \{ \Mb{h} \} } \subset \mathcal{F}_{[ \Mb{h}, \Mb{h}+1]}^{Z}.
\Eeq
Let $S \subset \mathds{Z}^d$. From \eqref{Eq_Inclusion_SigmafieldXk_SigmafieldZ}, it follows that
\Beq
\label{Eq_Inclusion_0_Sigmafield_Xtilde}
\sigma \left( \bigcup_{\Mb{h} \in S} \mathcal{F}^{\tilde{X}}_{ \{ \Mb{h} \} } \right) \subset \sigma \left( \bigcup_{\Mb{h} \in S} \mathcal{F}_{[ \Mb{h}, \Mb{h}+1]}^{Z} \right).
\Eeq
Additionally, it is easily shown that 
$$
\sigma \left( \bigcup_{\Mb{h} \in S} \mathcal{F}^{\tilde{X}}_{ \{ \Mb{h} \} } \right) = \mathcal{F}^{\tilde{X}}_{ S } \quad \mbox{and} \quad \sigma \left( \bigcup_{\Mb{h} \in S} \mathcal{F}^{Z}_{ [ \Mb{h}, \Mb{h}+1] } \right) = \mathcal{F}^{Z}_{ \bigcup_{\Mb{h} \in S} [ \Mb{h}, \Mb{h}+1]},
$$
which yield, using \eqref{Eq_Inclusion_0_Sigmafield_Xtilde}, that
\Beq
\label{Eq_Majoration_Sigmafield_Xtilde_Sigmafield_X}
\mathcal{F}^{\tilde{X}}_{ S } \subset \mathcal{F}^{Z}_{ \bigcup_{\Mb{h} \in S} [ \Mb{h}, \Mb{h}+1] }.
\Eeq
Thus, using \eqref{Eq_Def_2_Beta_Mixing} and \eqref{Eq_Majoration_Sigmafield_Xtilde_Sigmafield_X}, we obtain that, for all $S_1, S_2$ disjoint subsets of $\mathds{Z}^d$,
\Beq
\label{Eq_Majoration_Beta_Xtilde_Beta_Z}
\beta^{\tilde{X}} \left( S_1, S_2 \right) \leq \beta^Z \left( \bigcup_{\Mb{h}_1 \in S_1} [ \Mb{h}_1, \Mb{h}_1+1], \bigcup_{\Mb{h}_2 \in S_2} [ \Mb{h}_2, \Mb{h}_2+1] \right).
\Eeq

Now, $([ \Mb{h}_1, \Mb{h}_1+1])_{ \Mb{h}_1 \in S_1}$ and $([ \Mb{h}_2, \Mb{h}_2+1])_{ \Mb{h}_2 \in S_1}$ are countable families of compact subsets of $\mathds{R}^d$. Therefore, as $Z$ is a simple and sample-continuous max-stable random field on $\mathds{R}^d$, we can apply Theorem 2.2 in \cite{dombry2012strong}. The first point gives that, for any compact $S \subset \mathds{R}^d$, the quantity $C^Z(S)=\mathds{E} \left[ \sup_{\Mb{x} \in S} \left \{ Z(\Mb{x})^{-1} \right \} \right]$ is finite. Moreover, the third point yields that, for all $S_1, S_2$ disjoint subsets of $\mathds{Z}^d$,
\begin{align}
& \quad \ \beta^Z \left( \bigcup_{\Mb{h}_1 \in S_1} [ \Mb{h}_1, \Mb{h}_1+1], \bigcup_{\Mb{h}_2 \in S_2} [ \Mb{h}_2, \Mb{h}_2+1] \right) \nonumber
\\& \leq 2 \sum_{\Mb{h}_1 \in S_1} \sum_{\Mb{h}_2 \in S_2} \left[ C^Z \left( [ \Mb{h}_1, \Mb{h}_1+1] \right)+C^Z \left( [\Mb{h}_2, \Mb{h}_2+1] \right) \right] \nonumber
\\& \quad \quad \quad \quad \quad \quad \ \  \left[\theta^Z\left([ \Mb{h}_1, \Mb{h}_1+1] \right)+\theta^Z\left([\Mb{h}_2,\Mb{h}_2+1] \right)
-\theta^Z\left([ \Mb{h}_1, \Mb{h}_1+1]\cup[\Mb{h}_2,\Mb{h}_2+1] \right) \right].
\label{Eq_Majoration_Beta_Z_With_Extremal_Coefficients}
\end{align}
Let us introduce $K_2=C^Z \left( [0,1]^d \right)$. By stationarity of $Z$, we have that, for all $\Mb{h} \in \mathds{Z}^d$, 
\Beq
\label{Eq_Value_C_Cube_k}
C^Z \left( [\Mb{h}, \Mb{h}+1] \right)=K_2.
\Eeq
Furthermore, let $Y$ be a spectral random field of $Z$. Using the stationarity of $Z$, \eqref{Eq_Property_Extremal_Coefficient}, the linearity of the expectation and the fact that, for all $a, b \in \mathds{R}, a+b-\max \{a,b\}=\min\{a,b\}$, we have that
\begin{align}
& \quad \ \theta^Z\left([ \Mb{h}_1, \Mb{h}_1+1] \right)+\theta^Z\left([\Mb{h}_2,\Mb{h}_2+1] \right)
-\theta^Z\left([\Mb{h}_1, \Mb{h}_1+1] \cup[\Mb{h}_2,\Mb{h}_2+1] \right) \nonumber
\\& = \theta^Z\left([ 0,1]^d \right)+\theta^Z\left([\Mb{h}_2-\Mb{h}_1,\Mb{h}_2-\Mb{h}_1+1] \right)
-\theta^Z\left([ 0,1 ]^d \cup[\Mb{h}_2-\Mb{h}_1,\Mb{h}_2-\Mb{h}_1+1] \right) \nonumber
\\& = \mathds{E} \left[ \sup_{\Mb{x} \in [ 0,1]^d} \{ Y(\Mb{x}) \} + \sup_{\Mb{x} \in [\Mb{h}_2-\Mb{h}_1,\Mb{h}_2-\Mb{h}_1+1]} \{ Y(\Mb{x}) \} - \sup_{\Mb{x} \in [ 0,1]^d \bigcup [\Mb{h}_2-\Mb{h}_1,\Mb{h}_2-\Mb{h}_1 +1]} \{ Y(\Mb{x}) \}   \right] \nonumber
\\& = \mathds{E} \left[ \sup_{\Mb{x} \in [ 0,1]^d} \{ Y(\Mb{x}) \} + \sup_{\Mb{x} \in [\Mb{h}_2-\Mb{h}_1,\Mb{h}_2-\Mb{h}_1+1]} \{ Y(\Mb{x}) \} - \max \left \{ \sup_{\Mb{x} \in [ 0,1]^d} \{ Y(\Mb{x}) \},  \sup_{\Mb{x} \in [\Mb{h}_2-\Mb{h}_1,\Mb{h}_2-\Mb{h}_1+1]} \{ Y(\Mb{x}) \} \right \} \right] \nonumber
\\& = \mathds{E} \left[ \min \left \{ \sup_{\Mb{x} \in [ 0,1]^d} \{ Y(\Mb{x}) \},  \sup_{\Mb{x} \in [\Mb{h}_2-\Mb{h}_1,\Mb{h}_2-\Mb{h}_1+1]}  \{ Y(\Mb{x}) \} \right \} \right].
\label{Eq_Sum_Extremal_Coefficients}
\end{align}
Finally, combining \eqref{Eq_Majoration_Beta_Xtilde_Beta_Z}, \eqref{Eq_Majoration_Beta_Z_With_Extremal_Coefficients}, \eqref{Eq_Value_C_Cube_k} and \eqref{Eq_Sum_Extremal_Coefficients}, we obtain, for all $S_1, S_2$ disjoint subsets of $\mathds{Z}^d$, that
\Beq
\label{Eq_Majoration_Beta_Xtilde_With_Extremal_Coefficients}
\beta^{\tilde{X}} \left( S_1, S_2 \right) \leq 4K_2 \sum_{\Mb{h}_1 \in S_1} \sum_{\Mb{h}_2 \in S_2} \mathds{E} \left[ \min \left \{ \sup_{\Mb{x} \in [ 0,1]^d} \{ Y(\Mb{x}) \},  \sup_{\Mb{x} \in [\Mb{h}_2-\Mb{h}_1,\Mb{h}_2-\Mb{h}_1+1]}  \{ Y(\Mb{x}) \} \right \} \right].
\Eeq

Now, we show that \eqref{Eq_TCL2}, \eqref{Eq_TCL1} and \eqref{Eq_TCL3} are satisfied.
Using \eqref{Eq_Majoration_Alphamixing_With_Betamixing} and \eqref{eq:amix}, we obtain, for all $m\geq 1$ and $k,l\in\mathds{N}\cup\{\infty\}$, that 
\Beq
\label{Eq_Majoration_Alphakl_Xtilde_BetaklXtilde}
\alpha^{\tilde{X}}_{kl}(m) \leq \frac{1}{2} \sup \Big\{\beta^{\tilde{X}}(S_1,S_2):\ S_1, S_2 \subset \mathds{Z}^d,\ \ |S_1| \leq k, |S_2| \leq l,\ d(S_1,S_2)\geq m \Big\}.
\Eeq

Let $S_1$ and $S_2$ be subsets of $\mathds{Z}^d$ such that $|S_1| \leq k$, $|S_2| \leq l$ and $d(S_1, S_2) \geq m$, where $k, l \in \mathds{N}$ and $m \geq 1$. We have that
$$ \beta^{\tilde{X}}(S_1,S_2) \leq 4K_2kl \max_{\Mb{h}_1 \in S_1, \Mb{h}_2 \in S_2} \left \{ \mathds{E} \left[ \min \left \{ \sup_{\Mb{x} \in [0,1]^d} \{ Y(\Mb{x}) \},  \sup_{\Mb{x} \in [\Mb{h}_2-\Mb{h}_1, \Mb{h}_2-\Mb{h}_1+1]}  \{ Y(\Mb{x}) \} \right \} \right] \right \}.$$
Since $d(S_1, S_2) \geq m$, we have, for all $\Mb{h}_1 \in S_1$ and $\Mb{h}_2 \in S_2$, that $\| \Mb{h}_2-\Mb{h}_1 \| \geq m$. Thus, using \eqref{Eq_Decay_Cubes_Y}, we obtain that 
$$\beta^{\tilde{X}}(S_1,S_2) \leq 4KK_2kl m^{-b}.$$
Therefore, using \eqref{Eq_Majoration_Alphakl_Xtilde_BetaklXtilde}, for all $m \geq 1$ and $k, l \in \mathds{N}$, we have that
\Beq
\label{Eq_Majoration_alphakl_Xtilde}
\alpha^{\tilde{X}}_{kl}(m) \leq 2KK_2kl m^{-b}.
\Eeq
Hence, since $b > 2d$, we obtain, for all $k, l \geq 1$, that
$$ m^{d-1} \alpha^{\tilde{X}}_{kl}(m) \leq 2KK_2klm^{-(1+d)},$$
which immediately gives, since $d>0$ and $\alpha$-mixing coefficients are non-negative, that
$$ \sum_{m=1}^{\infty} m^{d-1} \alpha^{\tilde{X}}_{kl}(m) < \infty.$$
Thus, \eqref{Eq_TCL2} is satisfied.

Let $S_1$ and $S_2$ be subsets of $\mathds{Z}^d$ such that $|S_1| \leq k$, $|S_2| \leq l$ and $d(S_1, S_2) \geq m$, where $k \in \mathds{N}$, $l= \infty$ and $m \geq 1$. As $Y$ is positive, it is clear that, for all $\Mb{h}_1, \Mb{h}_2 \in \mathds{Z}^d$, 
$$ \mathds{E} \left[ \min \left \{ \sup_{\Mb{x} \in [ 0,1]^d} \{ Y(\Mb{x}) \},  \sup_{\Mb{x} \in [\Mb{h}_2-\Mb{h}_1,\Mb{h}_2-\Mb{h}_1+1]}  \{ Y(\Mb{x}) \} \right \} \right] \geq 0.$$
Hence, \eqref{Eq_Majoration_Beta_Xtilde_With_Extremal_Coefficients} gives that
$$ \beta^{\tilde{X}} \left( S_1, S_2 \right) \leq 4K_2 \sum_{\Mb{h}_1 \in S_1} \sum_{\Mb{h}_2 \in \mathds{Z}^d: \| \Mb{h}_2-\Mb{h}_1 \| \geq m} \mathds{E} \left[ \min \left \{ \sup_{\Mb{x} \in [ 0,1]^d} \{ Y(\Mb{x}) \},  \sup_{\Mb{x} \in [\Mb{h}_2-\Mb{h}_1,\Mb{h}_2-\Mb{h}_1+1]}  \{ Y(\Mb{x}) \} \right \} \right].$$
If follows from \eqref{Eq_Decay_Cubes_Y} that
$$ \beta^{\tilde{X}} \left( S_1, S_2 \right) \leq 4KK_2k \sum_{\Mb{h} \in \mathds{Z}^d: \| \Mb{h} \| \geq m} \| \Mb{h} \|^{-b}.$$
Consequently, \eqref{Eq_Majoration_Alphakl_Xtilde_BetaklXtilde} gives that
$$ \alpha_{1, \infty}^{\tilde{X}}(m) \leq 2KK_2k \sum_{\Mb{h} \in \mathds{Z}^d: \| \Mb{h} \| \geq m} \| \Mb{h} \|^{-b}.$$
Since $b>2d$ and $\alpha$-mixing coefficients are non-negative, we easily obtain that $\alpha_{1, \infty}^{\tilde{X}}(m)=o\left( m^{-d} \right)$ and hence that \eqref{Eq_TCL1} is satisfied.

Using H\"older's inequality and the fact that $\lambda([0,1]^d)=1$, we have that
\begin{align}
\left| \int_{[0,1]^d} X(\Mb{x}) \ \lambda(\mathrm{d}\Mb{x}) \right | &\leq \left( \int_{[ 0,1]^d} |X(\Mb{x})|^{2+\delta} \ \lambda(\mathrm{d}\Mb{x}) \right)^{\frac{1}{2+\delta}} \left( \int_{[ 0,1]^d} 1^{\frac{2+\delta}{1+\delta}} \ \lambda(\mathrm{d}\Mb{x}) \right)^{\frac{1+\delta}{2+\delta}} \nonumber
\\& = \left( \int_{[ 0,1]^d} |X(\Mb{x})|^{2+\delta} \ \lambda(\mathrm{d}\Mb{x}) \right)^{\frac{1}{2+\delta}}.
\nonumber
\end{align}
Hence, using \eqref{Eq_Def_Xtilde} and taking advantage of the stationarity of $X$, we have that
\begin{align*}
\mathds{E} \left[ |\tilde{X}(\Mb{0})|^{2+\delta} \right] \leq \mathds{E}\left[ \int_{[ 0,1]^d} |X(\Mb{x})|^{2+\delta} \ \lambda(\mathrm{d}\Mb{x}) \right] &=\int_{[0,1]^d} \mathds{E} \left[ |X(\Mb{x})|^{2+\delta} \right] \ \lambda(\mathrm{d}\Mb{x})
\\& =\mathds{E} \left[ |X(\Mb{0})|^{2+\delta} \right].
\end{align*}
Thus, using \eqref{Eq_Assumption_F}, we obtain that $\mathds{E}\left[ |\tilde{X}(\Mb{0})|^{2+\delta} \right]<\infty$.
Using \eqref{Eq_Majoration_alphakl_Xtilde}, we obtain that
$$ m^{d-1} \left( \alpha^{\tilde{X}}_{11}(m) \right)^{\frac{\delta}{2+\delta}} \leq m^{d-1} (2KK_2)^{\frac{\delta}{2+\delta}} m^{-\frac{b \delta }{2+\delta}}
= (2KK_2)^{\frac{\delta}{2+\delta}} m^{d-1-\frac{b \delta }{2+\delta}}.$$
Since $b > d(2+\delta)/\delta$, we have that $d-1-b \delta/(2+\delta)<-1$.
Therefore, since $\alpha$-mixing coefficients are non-negative, this yields
$$\sum_{m=1}^{\infty} m^{d-1} \left( \alpha^{\tilde{X}}_{11}(m) \right)^{\frac{\delta}{2+\delta}}< \infty.$$
Hence, \eqref{Eq_TCL3} is satisfied.

\medskip

Now, we recall that $\Lambda_n=\{ \Mb{h} \in \mathds{Z}^d: [\Mb{h}, \Mb{h}+1] \subset A_n \}$, where 
$$A_n= \bigcup_{\Mb{h} \in \mathds{Z}^d:[\Mb{h}, \Mb{h}+1] \subset V_n} [\Mb{h}, \Mb{h}+1].$$ Since, for any $n \in \mathds{N}$, $A_n$ is a bounded subset of $\mathds{R}^d$ (as a subset of $V_n$ which is bounded), and, by definition, $\Lambda_n \subset A_n$, we have that $\Lambda_n$ is a finite subset of $\mathds{Z}^d$. Moreover, as $(V_n)_{n \in \mathds{N}}$ is a Van Hove sequence, we have, for all $n \in \mathds{N}$, that $V_n \subset V_{n+1}$. This directly implies that $A_n \subset A_{n+1}$ and thus that $\Lambda_n \subset \Lambda_{n+1}$. Furthermore, by definition, we know that $\lim_{n \to \infty} V_n= \mathds{R}^d$. Hence, since, for all $n \in \mathds{N}$ and $\Mb{x} \in \partial A_n$, $\mathrm{dist}(\Mb{x},\partial V_n) \leq \sqrt{d}$, we obtain that $\lim_{n \to \infty} A_n= \mathds{R}^d$. Therefore, $\lim_{n \to \infty} \Lambda_n=\mathds{Z}^d$. Hence, $(\Lambda_n)_{n \in \mathds{N}}$ is a sequence of finite subsets of $\mathds{Z}^d$ which increases to $\mathds{Z}^d$.
Furthermore, we have that, for all $n \in \mathds{N}$, 
\Beq
\label{Eq_Inclusion_Remainder_Tube}
V_n \backslash A_n \subset N \left( \partial V_n, \sqrt{d} \right).
\Eeq
Hence, $\lambda(V_n \backslash A_n) \leq \lambda \left( N \left( \partial V_n, \sqrt{d} \right) \right)$, which gives, since $(V_n)_{n \in \mathds{N}}$ is a Van Hove sequence, that
\Beq
\label{Eq_Lim_LambdaReminder_On_LambdaVn}
\lim_{n \to \infty} \frac{\lambda(V_n \backslash A_n)}{\lambda(V_n)}=0.
\Eeq
Since $\lambda(A_n)=\lambda(V_n)-\lambda(V_n \backslash A_n)$, we have that
$$ \frac{\lambda(A_n)}{\lambda(V_n)}=1-\frac{\lambda(V_n \backslash A_n)}{\lambda(V_n)}.$$ 
Thus, using \eqref{Eq_Lim_LambdaReminder_On_LambdaVn}, we obtain that 
\Beq
\label{Eq_Lim_LambdaAn_On_LambdaVn}
\lim_{n \to \infty} \frac{\lambda(A_n)}{\lambda(V_n)}=1.
\Eeq
Moreover, since, for all $\Mb{h} \in \mathds{Z}^d$, $\lambda([\Mb{h}, \Mb{h}+1])=1$, we obtain that, for all $n \in \mathds{N}$, 
\Beq
\label{Eq_Equality_Cardinal_Lambdan_Measure_An}
|\Lambda_n|=\lambda(A_n).
\Eeq
For the same reason, we have, for all $n \in \mathds{N}$, that $|\partial \Lambda_n| \leq \lambda \left( N \left( \partial A_n, \sqrt{d} \right) \right)$. Additionally, as $\mathrm{dist}(\partial V_n, \partial A_n) \leq \sqrt{d}$, we obtain that $N \left( \partial A_n, \sqrt{d} \right) \subset N \left( \partial V_n, 2\sqrt{d} \right)$. Hence, we have, for all $n \in \mathds{N}$, that $|\partial \Lambda_n| \leq \lambda \left( N \left( \partial V_n, 2\sqrt{d} \right) \right)$.
Therefore, using \eqref{Eq_Equality_Cardinal_Lambdan_Measure_An}, it follows that
$$ \frac{| \partial \Lambda_n |}{|\Lambda_n|}=\frac{|\partial \Lambda_n|}{\lambda(V_n)} \frac{\lambda(V_n)}{\lambda(A_n)} \leq \frac{\lambda \left( N \left( \partial V_n, 2\sqrt{d} \right) \right)}{\lambda(V_n)} \frac{\lambda(V_n)}{\lambda(A_n)}.$$
Consequently, using \eqref{Eq_Lim_LambdaAn_On_LambdaVn} and the fact that $(V_n)_{n \in \mathds{N}}$ is a Van Hove sequence, we obtain that $\lim_{n \to \infty} |\partial \Lambda_n|/|\Lambda_n|=0$. To summarize, we have that $(\Lambda_n)_{n \in \mathds{N}}$ is a sequence of finite subsets of $\mathds{Z}^d$ which increases to $\mathds{Z}^d$ and such that $\lim_{n \to \infty} |\partial \Lambda_n|/|\Lambda_n|=0$. 

\medskip

Thus, Theorem \ref{theo:BCLT} gives that $\sum_{\Mb{h} \in \mathds{Z}^d} \left | \mathrm{Cov} \left(\tilde{X}(\Mb{0}), \tilde{X}(\Mb{h})\right) \right | < \infty$. We introduce $\sigma_1^2=\sum_{\Mb{h} \in \mathds{Z}^d} \mathrm{Cov} \left(\tilde{X}(\Mb{0}), \tilde{X}(\Mb{h})\right)$. 
Using the fact that the covariance is bilinear, the stationarity of $X$ and the definition of $\sigma^2$ in \eqref{Eq_Def_Sigma_Square}, we have that
\begin{align}
\sigma_1^2 & = \int_{[0,1]^d} \left( \sum_{\Mb{h} \in \mathds{Z}^d}  \int_{[\Mb{h}, \Mb{h}+1]} \mathrm{Cov}(X(\Mb{x}), X(\Mb{y})) \ \lambda(\mathrm{d}\Mb{y}) \right) \lambda(\mathrm{d}\Mb{x}) \nonumber
\\& = \int_{[0,1]^d} \left( \int_{\mathds{R}^d} \mathrm{Cov}(X(\Mb{x}), X(\Mb{y})) \ \lambda(\mathrm{d}\Mb{y}) \right) \ \lambda(\mathrm{d}\Mb{x}) \nonumber
\\& = \int_{[0,1]^d} \left( \int_{\mathds{R}^d} \mathrm{Cov}(X(\Mb{0}), X(\Mb{y}-\Mb{x})) \ \lambda(\mathrm{d}\Mb{y}) \right) \ \lambda(\mathrm{d}\Mb{x}) \nonumber
\\& = \sigma^2.
\label{Eq_Equality_Covariance_Xtilde_Covariance}
\end{align}
Consequently, it follows from \eqref{Eq_Condition_Positivity_Sigma2} that $\sigma_1^2>0$. Hence, Theorem \ref{theo:BCLT} yields that
\Beq
\label{Eq_TCL_Xtilde}
\frac{1}{|\Lambda_n|^{\frac{1}{2}}} \sum_{\Mb{h} \in \Lambda_n} \tilde{X}(\Mb{h}) \overset{d}{\rightarrow} \mathcal{N}(0,\sigma_1^2). 
\Eeq
Finally, combining \eqref{Eq_Def_In1_In2}, \eqref{Eq_Integral_Sum_Representation}, \eqref{Eq_Equality_Cardinal_Lambdan_Measure_An}, \eqref{Eq_Equality_Covariance_Xtilde_Covariance} and \eqref{Eq_TCL_Xtilde}, we obtain that
$$ \left( \frac{\lambda(V_n)}{\lambda(A_n)} \right)^{\frac{1}{2}} I_{n,1} \overset{d}{\rightarrow} \mathcal{N}(0, \sigma^2), \mbox{ for } n \to \infty.$$
Hence, at last, using \eqref{Eq_Lim_LambdaAn_On_LambdaVn}, Slutsky's theorem yields that
\Beq
\label{Eq_Convergence_In1}
I_{n,1} \overset{d}{\rightarrow} \mathcal{N}(0, \sigma^2), \mbox{ for } n \to \infty.
\Eeq 

\medskip
\medskip

\noindent \textbf{Part 3: Study of $(I_{n,2})_{n \in \mathds{N}}$}

\medskip
\medskip

We now focus on the second term in \eqref{Eq_Decomposition_Integral_Total}, i.e. $I_{n,2}$. Using \eqref{Eq_Inclusion_Remainder_Tube}, the stationarity of $X$ and \eqref{Eq_Convergence_Int_Abs_Cov}, we obtain that, for all $n \in \mathds{N}$,
\begin{align}
\mathrm{Var}(I_{n,2})&=\frac{1}{\lambda(V_n)} \mathrm{Var} \left( \int_{V_n \backslash A_n} X(\Mb{x}) \ \lambda(\mathrm{d}\Mb{x}) \right) \nonumber
\\& = \frac{1}{\lambda(V_n)}  \int_{V_n \backslash A_n}  \int_{V_n \backslash A_n} \mathrm{Cov}(X(\Mb{x}), X(\Mb{y}))  \ \lambda(\mathrm{d}\Mb{x})  \ 
\lambda(\mathrm{d}\Mb{y}) \nonumber
\\& \leq \frac{1}{\lambda(V_n)}  \int_{N \left( \partial V_n, \sqrt{d} \right)}  \int_{N \left( \partial V_n, \sqrt{d} \right)} | \mathrm{Cov}(X(\Mb{x}), X(\Mb{y}))| \ \lambda(\mathrm{d}\Mb{x})  \ \lambda(\mathrm{d}\Mb{y}) \nonumber
\\& = \frac{1}{\lambda(V_n)}  \int_{N \left( \partial V_n, \sqrt{d} \right)} \left( \int_{N \left( \partial V_n, \sqrt{d} \right)} | \mathrm{Cov}(X(\Mb{0}), X(\Mb{y}-\Mb{x})) |  \ \lambda(\mathrm{d}\Mb{y}) \right)  \ \lambda(\mathrm{d}\Mb{x}) \nonumber
\\& \leq \frac{\lambda \left( N \left( \partial V_n, \sqrt{d} \right) \right)}{\lambda(V_n)} K_1 \nonumber.
\end{align}
Since $(V_n)_{n \in \mathds{N}}$ is a Van Hove sequence, we have that $\lim_{n \to \infty} \lambda \left( N \left( \partial V_n, \sqrt{d} \right) \right)/\lambda(V_n)=0$, giving that $\lim_{n \to \infty} \mathrm{Var}(I_{n,2})=0$. Since $\mathds{E}[I_{n,2}]=0$, this shows (using Bienaym\'e-Tchebychev's inequality) that $(I_{n,2})_{n \in \mathds{N}}$ tends towards $0$ in probability. 

\medskip
\medskip

Finally, using \eqref{Eq_Decomposition_Integral_Total} and \eqref{Eq_Convergence_In1} and applying Slutsky's theorem, we obtain that
$$\frac{1}{[\lambda(V_n)]^{\frac{1}{2}}} \int_{V_n} X(\Mb{x}) \ \lambda(\mathrm{d}\Mb{x}) \overset{d}{\rightarrow} \mathcal{N}(0, \sigma^2), \mbox{ for } n \to \infty.$$
This completes the proof.
\end{proof}

\subsection{For Proposition \ref{Prop_Condition_Sigma_Positive}}

\begin{proof}
Since the random field $Z$ is max-stable, we know that, for all $\Mb{x} \in \mathds{R}^d$, the random vector $\Mb{Z}=(Z(\Mb{0}), Z(\Mb{x}))^{'}$ is max-stable and thus max-infinitely divisible. Hence, Proposition 5.29 in \cite{resnickextreme} gives that it is associated. 
In Definition \ref{Def_Association}, let us choose $q=2$ and define $g_1$ and $g_2$ as follows:
$$ g_1(z_1,z_2)=F(z_1) \quad \mathrm{and} \quad g_2(z_1,z_2)=F(z_2), \quad z_1, z_2 \in \mathds{R}.$$ 
As $F$ is non-decreasing, $g_1$ and $g_2$ are non-decreasing in the sense of Definition \ref{Def_Association}.
Moreover, since $F$ satisfies \eqref{Eq_Assumption_F}, we have that $E[|F(Z(\Mb{0}))|^2] < \infty$, $E[|F(Z(\Mb{x}))|^2] < \infty$ and, using Cauchy-Schwarz inequality, $E[|F(Z(\Mb{0}))F(Z(\Mb{x}))|] < \infty$. This implies that $\mathds{E}[|g_1(\Mb{Z})|^2] < \infty$, $\mathds{E}[|g_2(\Mb{Z})|^2] < \infty$ and $\mathds{E}[|g_1(\Mb{Z}) g_2(\Mb{Z})|] < \infty$. By definition of association, it follows that, for all $\Mb{x} \in \mathds{R}^d$, $\mathrm{Cov}(F(Z(\Mb{0}))), F(Z(\Mb{x}))) \geq 0$, i.e.
\Beq
\label{Eq_Nonnegativity_Covariance}
\mathrm{Cov}(X(\Mb{0}), X(\Mb{x})) \geq 0.
\Eeq

Now, since $Z$ is max-stable and $F$ is measurable, non-decreasing and non-constant, we have that 
\Beq
\label{Eq_Positivity_Var}
\mathrm{Var}(X(\Mb{0}))>0.
\Eeq

Since $F$ is monotone, the set of points at which $F$ is not continuous, denoted $\mathcal{D}_F$, is at most countable. Hence, for all $\Mb{x}_0 \in \mathds{R}^d$, since $Z(\Mb{x}_0)$ is a continuous random variable (standard Fr\'echet), we have that $\mathds{P}(Z(\Mb{x}_0) \in \mathcal{D}_F )=0$. Thus, as $Z$ is sample-continuous, $X$ is almost surely (a.s.) continuous at $\Mb{x}_0$, which implies that, for all $\Mb{x}_0 \in \mathds{R}^d$, 
\Beq
\label{Eq_Convergence_Square_Diff_To_Zero}
\mathds{P} \left( \lim_{\Mb{x} \to \Mb{x}_0} |X(\Mb{x})-X(\Mb{x}_0)|^2=0 \right)=1.
\Eeq
We introduce $\delta_1=\delta/2$, where $\delta$ appears in \eqref{Eq_Assumption_F}. Using the well-known fact that, for all $a,b \in \mathds{R}$ and $p \geq 1$, $|a-b|^p \leq 2^{p-1} (|a|^p+|b|^p)$, we obtain that 
$$ \left (|X(\Mb{x})-X(\Mb{x}_0)|^2 \right)^{1+\delta_1} \leq 2^{1+\delta} (|X(\Mb{x})|^{2+\delta} + |X(\Mb{x}_0)|^{2+\delta} ).$$
Using the stationarity of $X$ and \eqref{Eq_Assumption_F}, we obtain, for all $\Mb{x}_0 \in \mathds{R}^d$, that 
$$ \sup_{\Mb{x} \in \mathds{R}^d} \left \{ \mathds{E} \left[ \left (|X(\Mb{x})-X(\Mb{x}_0)|^2 \right)^{1+\delta_1} \right] \right \} \leq 2^{2+\delta} \mathds{E}\left[ X(\Mb{0}) \right]^{2+\delta}<\infty,$$
implying, since $\delta_1>0$, that the random field $\left \{ |X(\Mb{x})-X(\Mb{x}_0)|^2 \right \}_{\Mb{x} \in \mathds{R}^d}$ is uniformly integrable.
Consequently, we obtain using \eqref{Eq_Convergence_Square_Diff_To_Zero} that, for all $\Mb{x}_0 \in \mathds{R}^d$, $\lim_{\Mb{x} \to \Mb{x}_0} \mathds{E}[|X(\Mb{x})-X(\Mb{x}_0)|^2]=0$, meaning that $X$ is continuous in quadratic mean. Hence, since $X$ is second-order stationary (since it is strictly stationary and has a second moment), its covariance function, defined by $\mathrm{Cov}_X(\Mb{x})=\mathrm{Cov}(X(\Mb{0}), X(\Mb{x}))$, $\Mb{x} \in \mathds{R}^d$, is continuous at the origin. Hence, there exists $\xi>0$ such that, for all $\Mb{x} \in \mathds{R}^d$ satisfying $\| \Mb{x} \| \leq \xi$, $|\mathrm{Cov}(X(\Mb{0}), X(\Mb{x}))-\mathrm{Var}(X(\Mb{0}))| \leq \mathrm{Var}(X(\Mb{0}))/2$, implying that 
$\mathrm{Cov}(X(\Mb{0}), X(\Mb{x})) \geq \mathrm{Var}(X(\Mb{0}))/2$. Thus, using \eqref{Eq_Nonnegativity_Covariance}, \eqref{Eq_Positivity_Var} and the fact that $\xi>0$, we obtain that
$$\sigma^2 \geq \int_{\Mb{x} \in \mathds{R}^d: \| \Mb{x} \| \leq \xi} \mathrm{Cov}(X(\Mb{0}), X(\Mb{x})) \ \lambda (\mathrm{d}\Mb{x}) \geq \lambda \left( \left \{\Mb{x} \in \mathds{R}^d: \| \Mb{x} \| \leq \xi \right \} \right) \frac{\mathrm{Var}(X(\Mb{0}))}{2} > 0.$$
This concludes the proof.
\end{proof}

\subsection{For Proposition \ref{Prop_Upper_Bound_Expectation_Min_Brown_Resnick}}

\begin{proof}
We introduce 
\Beq
\label{Eq_Definition_Eta}
\nu(\Mb{h})=\mathds{E} \left[ \min \left \{ \sup_{\Mb{x} \in [0,1]^d} \{ Y(\Mb{x}) \},  \sup_{\Mb{x} \in [\Mb{h},\Mb{h}+1]}  \{ Y(\Mb{x}) \} \right \} \right], \quad \Mb{h} \in \mathds{Z}^d.
\Eeq
Since $Y$ is positive, we have, for all $\Mb{h} \in \mathds{Z}^d$, that
\Beq
\label{Eq_Positivity_Nu}
\nu(\Mb{h}) \geq 0.
\Eeq
Moreover, we have, for all $\Mb{h} \in \mathds{Z}^d$, that
\begin{align}
\nu(\Mb{h}) &=\mathds{E} \left[ \min \left \{ \sup_{\Mb{x} \in [0,1]^d} \{ Y(\Mb{x}) \},  \sup_{\Mb{x} \in [0,1]^d}  \{ Y(\Mb{x}+\Mb{h}) \} \right \} \right] \nonumber
\\& = \mathds{E} \left[ \min \left \{ \sup_{\Mb{x} \in [0,1]^d} \{ Y(\Mb{x}) \}, Y(\Mb{h}) \sup_{\Mb{x} \in [0,1]^d}  \left \{ \frac{Y(\Mb{x}+\Mb{h})}{Y(\Mb{h})} \right \} \right \} \right].
\label{Eq_First_Computation_Nu}
\end{align}
Let $\varepsilon$ denote any function from $\mathds{Z}^d$ to $(0,\infty)$. 
We have that
\begin{align*}
& \quad \ \sup_{\Mb{x} \in [0,1]^d} \{ Y(\Mb{x}) \} \ \mathds{I}_{ \{ Y(\Mb{h})>\varepsilon(\Mb{h}) \} } + \varepsilon(\Mb{h}) \sup_{\Mb{x} \in [0,1]^d}  \left \{ \frac{Y(\Mb{x}+\Mb{h})}{Y(\Mb{h})} \right \} \  \mathds{I}_{ \{ Y(\Mb{h}) \leq\varepsilon(\Mb{h}) \} }
\\& = \left \{ 
\begin{array}{cc}
\sup_{\Mb{x} \in [0,1]^d} \{ Y(\Mb{x}) \} & \mbox{if } Y(\Mb{h})>\varepsilon(\Mb{h}), \\
\varepsilon(\Mb{h}) \sup_{\Mb{x} \in [0,1]^d}  \left \{ \frac{Y(\Mb{x}+\Mb{h})}{Y(\Mb{h})} \right \} & \mbox{if } Y(\Mb{h}) \leq \varepsilon(\Mb{h}),
\end{array}
\right.
\end{align*}
yielding, since $Y$ is positive, that 
\begin{align}
& \quad \ \  \min \left \{ \sup_{\Mb{x} \in [0,1]^d} \{ Y(\Mb{x}) \}, Y(\Mb{h}) \sup_{\Mb{x} \in [0,1]^d}  \left \{ \frac{Y(\Mb{x}+\Mb{h})}{Y(\Mb{h})} \right \} \right \}  \nonumber
\\& \leq \sup_{\Mb{x} \in [0,1]^d} \{ Y(\Mb{x}) \} \ \mathds{I}_{ \{ Y(\Mb{h})>\varepsilon(\Mb{h}) \} } + \varepsilon(\Mb{h}) \sup_{\Mb{x} \in [0,1]^d}  \left \{ \frac{Y(\Mb{x}+\Mb{h})}{Y(\Mb{h})} \right \} \  \mathds{I}_{ \{ Y(\Mb{h}) \leq\varepsilon(\Mb{h}) \} }.
\label{Eq_Explanation_Majoration_Epsilon}
\end{align}
We obtain, using \eqref{Eq_First_Computation_Nu}, \eqref{Eq_Explanation_Majoration_Epsilon} and Cauchy-Schwarz inequality, that, for all $\Mb{h} \in \mathds{Z}^d$,
\begin{align}
\nu(\Mb{h}) &\leq \mathds{E} \left[ \sup_{\Mb{x} \in [0,1]^d} \{ Y(\Mb{x}) \} \ \mathds{I}_{ \{ Y(\Mb{h})>\varepsilon(\Mb{h}) \} } + \varepsilon(\Mb{h}) \sup_{\Mb{x} \in [0,1]^d}  \left \{ \frac{Y(\Mb{x}+\Mb{h})}{Y(\Mb{h})} \right \} \  \mathds{I}_{ \{ Y(\Mb{h}) \leq\varepsilon(\Mb{h}) \} } \right] \nonumber 
\\& \leq \mathds{E} \left[ \sup_{\Mb{x} \in [0,1]^d} \{ Y(\Mb{x}) \} \ \mathds{I}_{ \{ Y(\Mb{h})>\varepsilon(\Mb{h}) \} } \right] + \mathds{E} \left[ \varepsilon(\Mb{h}) \sup_{\Mb{x} \in [0,1]^d}  \left \{ \frac{Y(\Mb{x}+\Mb{h})}{Y(\Mb{h})} \right \} \  \mathds{I}_{ \{ Y(\Mb{h}) \leq\varepsilon(\Mb{h}) \} } \right] \nonumber
\\& \leq \mathds{E} \left [ \sup_{\Mb{x} \in [0,1]^d} \left \{ Y^2(\Mb{x}) \right \} \right ]^{\frac{1}{2}} \mathds{P}(Y(\Mb{h}) > \varepsilon(\Mb{h}) )^{\frac{1}{2}} + \varepsilon(\Mb{h}) \mathds{E} \left [ \sup_{\Mb{x} \in [0,1]^d}  \left \{ \frac{Y(\Mb{x}+\Mb{h})}{Y(\Mb{h})} \right \} \right].
\label{Eq_First_Majoration_Gamma}
\end{align}
Since $W$ is a Gaussian random field, we have that 
\Beq
\label{Eq_Gaussian_Survivor}
\mathds{P}(Y(\Mb{h}) > \varepsilon(\Mb{h})) = \mathds{P} \left( W(\Mb{h}) > \frac{\sigma_W^2(\Mb{h})}{2} + \log(\varepsilon(\Mb{h})) \right)=\bar{\Phi} \left( \frac{\sigma_W^2(\Mb{h})}{2} + \log(\varepsilon(\Mb{h})) \right),
\Eeq
where $\bar{\Phi}=1-\Phi$ with $\Phi$ being the standard Gaussian distribution function.
Now, for all $\Mb{h} \in \mathds{Z}^d$ and $\Mb{x} \in [0,1]^d$, we have that
$$ \frac{Y(\Mb{x}+\Mb{h})}{Y(\Mb{h})} = \exp \left( W(\Mb{x}+\Mb{h})-W(\Mb{h})-\frac{\sigma_W^2(\Mb{x}+\Mb{h})-\sigma_W^2(\Mb{h})}{2} \right).$$
Hence, since $W$ has stationary increments, we obtain that, for all $\Mb{h} \in \mathds{Z}^d$,
$$ \left \{ \frac{Y(\Mb{x}+\Mb{h})}{Y(\Mb{h})} \right \}_{\Mb{x} \in \mathds{R}^d} \overset{d}{=} \left \{ \exp \left( W(\Mb{x})-\frac{\sigma_W^2(\Mb{x}+\Mb{h})-\sigma_W^2(\Mb{h})}{2} \right) \right \}_{\Mb{x} \in \mathds{R}^d},$$
which yields
\Beq
\label{Eq_Equality_Distribution_Fraction}
\left \{ \frac{Y(\Mb{x}+\Mb{h})}{Y(\Mb{h})} \right \}_{\Mb{x} \in \mathds{R}^d} \overset{d}{=} \left \{ \exp \left( W(\Mb{x}) -\frac{\sigma_W^2(\Mb{x})}{2} \right) \exp \left( \frac{\sigma_W^2(\Mb{x}) + \sigma_W^2(\Mb{h})-\sigma_W^2(\Mb{x}+\Mb{h})}{2} \right) \right \}_{\Mb{x} \in \mathds{R}^d}.
\Eeq
Now, we show that $\mathds{E} \left[ \sup_{\Mb{x} \in [0,1]^d} \{ Y(\Mb{x}) \} \right] < \infty.$
Using the fact that the exponential function is increasing, we have that 
\Beq
\label{Eq_Expectation_Y}
\mathds{E} \left[ \sup_{\Mb{x} \in [0,1]^d} \{ Y(\Mb{x}) \} \right] = \mathds{E} \left[ \exp \left( \sup_{\Mb{x} \in [0,1]^d} \left \{ W(\Mb{x})-\frac{\sigma_W^2(\Mb{x})}{2} \right \} \right) \right] \leq \mathds{E} \left [ \exp \left( \sup_{\Mb{x} \in [0,1]^d} \left \{ W(\Mb{x}) \right \} \right) \right].
\Eeq
As $W$ is a centred Gaussian random field which is a.s. bounded on $[0,1]^d$, Theorem 2.1.2 in \cite{adler2007random} yields that 
$\sup_{\Mb{x} \in [0,1]^d} \left \{ \mathds{E} \left[ (W(\Mb{x}))^2 \right] \right \} < \infty.$
Hence, let us introduce
\Beq
\label{Eq_Def_Tau}
\tau= \sup_{\Mb{x} \in [0,1]^d} \left \{ \mathds{E} \left[ (W(\Mb{x}))^2 \right] \right \}.
\Eeq
It is clear that $\tau >0$. Since $W$ is a centred Gaussian random field which is a.s. bounded on $[0,1]^d$, Theorem 2.1.1 in \cite{adler2007random} gives that, for all $u>0$,
$$ \mathds{P} \left ( \sup_{\Mb{x} \in [0,1]^d} \left \{ W(\Mb{x}) \right \} > u \right ) \leq  \exp \left( -\frac{u^2}{2 \tau} \right),$$
which yields, for all $w >0$, that
\Beq
\label{Eq_Majoration_Survival_Exp_Sup_W}
\mathds{P} \left ( \exp \left( \sup_{\Mb{x} \in [0,1]^d} \left \{ W(\Mb{x}) \right \} \right) > w \right ) \leq  \exp \left( -\frac{\ln^2(w)}{2 \tau} \right).
\Eeq
Using \eqref{Eq_Majoration_Survival_Exp_Sup_W} and two changes of variable, we obtain that
\begin{align}
\mathds{E} \left [ \exp \left( \sup_{\Mb{x} \in [0,1]^d} \left \{ W(\Mb{x}) \right \} \right) \right] &=\int_{0}^{\infty} \mathds{P} \left ( \exp \left( \sup_{\Mb{x} \in [0,1]^d} \left \{ W(\Mb{x}) \right \} \right) > w \right ) \ \mathrm{d}w \nonumber
\\& \leq \int_{- \infty}^{\infty} \exp \left( -\frac{v^2}{2 \tau} + v \right) \  \mathrm{d}v  \nonumber
\\& = \exp \left( \frac{\tau}{2} \right) \int_{- \infty}^{\infty} \exp \left( -\frac{v_1^2}{2 \tau} \right) \ \mathrm{d}v_1 < \infty.
\label{Eq_Majoration_Expectation_Exp_Sup_W}
\end{align}
Combining \eqref{Eq_Expectation_Y} and \eqref{Eq_Majoration_Expectation_Exp_Sup_W}, we obtain that $\mathds{E} \left[ \sup_{\Mb{x} \in [0,1]^d} \{ Y(\Mb{x}) \} \right] < \infty$. Very similar arguments yield that $\mathds{E} \left[ \sup_{\Mb{x} \in [0,1]^d} \left \{ Y^2(\Mb{x}) \right \} \right] < \infty$. Hence, we introduce 
\Beq
\label{Def_C1_C2}
C_1 = \mathds{E} \left[ \sup_{\Mb{x} \in [0,1]^d} \{ Y(\Mb{x}) \} \right] \quad \mbox{and} \quad C_2 = \mathds{E} \left[ \sup_{\Mb{x} \in [0,1]^d} \left \{ Y^2(\Mb{x}) \right \} \right].
\Eeq
The random fields $Y$ and $Y^2$ being positive, we have $C_1, C_2>0$. Furthermore, let
\Beq
\label{Eq_Def_Delta}
\delta(\Mb{h})=\sup_{\Mb{x} \in [0,1]^d} \left \{ \sigma_W^2(\Mb{x}) + \sigma_W^2(\Mb{h})-\sigma_W^2(\Mb{x}+\Mb{h}) \right \}, \quad \Mb{h} \in \mathds{Z}^d.
\Eeq
The combination of \eqref{Eq_First_Majoration_Gamma}, \eqref{Eq_Gaussian_Survivor}, \eqref{Eq_Equality_Distribution_Fraction}, \eqref{Def_C1_C2} and \eqref{Eq_Def_Delta} gives that, for all $\Mb{h} \in \mathds{Z}^d$,
$$ \nu(\Mb{h}) \leq C_2^{\frac{1}{2}} \bar{\Phi}^{\frac{1}{2}} \left( \frac{\sigma_W^2(\Mb{h})}{2} + \log(\varepsilon(\Mb{h})) \right) + \varepsilon(\Mb{h}) C_1 \exp \left( \frac{\delta(\Mb{h})}{2} \right).$$
Let us take $\varepsilon(\Mb{h})=\exp \left( -\frac{\sigma_W^2(\Mb{h})}{4} \right), \Mb{h} \in \mathds{Z}^d$. Hence, we obtain
\Beq
\label{Eq_Second_Majoration_Gamma}
\nu(\Mb{h}) \leq C_2^{\frac{1}{2}} \bar{\Phi}^{\frac{1}{2}} \left( \frac{\sigma_W^2(\Mb{h})}{4} \right) +C_1 \exp \left( \frac{\delta(\Mb{h})}{2} -\frac{\sigma_W^2(\Mb{h})}{4} \right).
\Eeq

We obtain using \eqref{Eq_Condition1_Variogram_TCL} that
$$ \delta(\Mb{h}) \leq \sup_{\Mb{x} \in [0,1]^d} \{ \sigma_W^2(\Mb{x}) \} + \sup_{\Mb{x} \in [0,1]^d} \{ \sigma_W^2(\Mb{h})-\sigma_W^2(\Mb{x}+\Mb{h}) \} \leq \sup_{\Mb{x} \in [0,1]^d} \{ \sigma_W^2(\Mb{x}) \} + o(\sigma_W^2(\Mb{h})).$$
As $W$ is centred, it follows from \eqref{Eq_Def_Tau} that $\sup_{\Mb{x} \in [0,1]^d} \{ \sigma_W^2(\Mb{x}) \}=\tau$. Thus,
\Beq
\label{Eq_Inequality_Difference_Delta_Sigma2Over4}
\frac{\delta(\Mb{h})}{2} -\frac{\sigma_W^2(\Mb{h})}{4} \leq \frac{\tau}{2} + o \left(\sigma_W^2(\Mb{h})\right) - \frac{\sigma_W^2(\Mb{h})}{4}.
\Eeq
Since $\lim_{\| \Mb{h} \| \to \infty} \sigma_W^2(\Mb{h})=\infty$, it is clear that 
$$ \frac{\tau}{2} + o(\sigma_W^2(\Mb{h})) - \frac{\sigma_W^2(\Mb{h})}{4} \underset{\| \Mb{h} \| \to \infty}{\sim} - \frac{\sigma_W^2(\Mb{h})}{4},$$
and hence that there exist $A, A_1>0$ such that, for all $\Mb{h} \in \mathds{Z}^d$ satisfying $\| \Mb{h} \| \geq A$,
$$
\frac{\tau}{2} + o(\sigma_W^2(\Mb{h})) - \frac{\sigma_W^2(\Mb{h})}{4} \leq -A_1 \frac{\sigma_W^2(\Mb{h})}{4}.
$$
Therefore, using \eqref{Eq_Inequality_Difference_Delta_Sigma2Over4}, we obtain, for all $\Mb{h} \in \mathds{Z}^d$ satisfying $\| \Mb{h} \| \geq A$, that
\Beq
\label{Eq_Upper_Bound_Delta_Difference_Delta_Sigma2Over4}
\exp \left( \frac{\delta(\Mb{h})}{2} - \frac{\sigma_W^2(\Mb{h})}{4} \right) \leq \exp \left( -A_1 \frac{\sigma_W^2(\Mb{h})}{4} \right).
\Eeq

Now, Mill's ratio gives that
$$ \bar{\Phi} \left( \frac{\sigma_W^2(\Mb{h})}{4} \right) \underset{\| \Mb{h} \| \to \infty}{\sim} \frac{4 \exp \left( - \frac{\sigma_W^4(\Mb{h})}{32} \right)}{(2\pi)^{\frac{1}{2}} \sigma_W^2(\Mb{h})} \mbox{ and thus } \bar{\Phi}^{\frac{1}{2}} \left( \frac{\sigma_W^2(\Mb{h})}{4} \right) \underset{\| \Mb{h} \| \to \infty}{\sim} \frac{2 \exp \left( - \frac{\sigma_W^4(\Mb{h})}{64} \right)}{(2\pi)^{\frac{1}{4}} \sigma_W(\Mb{h})}.$$
Hence, we easily obtain that there exists $A_2>0$ such that, for all $\Mb{h} \in \mathds{Z}^d$,
\Beq
\label{Eq_Upper_Bound_Phi_Bar}
\bar{\Phi}^{\frac{1}{2}} \left( \frac{\sigma_W^2(\Mb{h})}{4} \right) \leq A_2 \frac{\exp \left( - \frac{\sigma_W^4(\Mb{h})}{64}\right)}{\sigma_W(\Mb{h})}.
\Eeq
Combining \eqref{Eq_Second_Majoration_Gamma}, \eqref{Eq_Upper_Bound_Delta_Difference_Delta_Sigma2Over4} and \eqref{Eq_Upper_Bound_Phi_Bar}, we obtain that there exists $A_3>0$ such that, for all $\Mb{h} \in \mathds{Z}^d$,
\Beq
\label{Eq_Second_Majoration_Nu}
\nu(\Mb{h}) \leq A_3 \exp \left( -A_1 \frac{\sigma_W^2(\Mb{h})}{4} \right).
\Eeq
Using \eqref{Eq_Condition2_Variogram_TCL}, we have, for all $b>0$, that
$$ \lim_{\| \Mb{h} \| \to \infty} \frac{\exp \left( -A_1 \frac{\sigma_W^2(\Mb{h})}{4} \right)}{\| \Mb{h} \|^{-b}}= \lim_{\| \Mb{h} \| \to \infty} \exp \left( -A_1 \frac{\sigma_W^2(\Mb{h})}{4}+b \ln(\| \Mb{h} \| ) \right)=0,$$
giving, using \eqref{Eq_Positivity_Nu} and \eqref{Eq_Second_Majoration_Nu}, that $\nu(\Mb{h}) \underset{\| \Mb{h} \| \to \infty}{=} o(\| \Mb{h} \|^{-b})$. Thus, for all $b>0$, there exists $A_4>0$ such that 
\Beq
\label{Eq_Final_Upper_Bound_Eta_Brown_Resnick}
\nu(\Mb{h}) \leq A_4 \| \Mb{h} \|^{-b}.
\Eeq
Finally, combining \eqref{Eq_Definition_Eta} and \eqref{Eq_Final_Upper_Bound_Eta_Brown_Resnick}, we obtain that
\eqref{Eq_Decay_Cubes_Y} is satisfied for all $b>0$.
\end{proof}

\subsection{For Proposition \ref{Prop_Variogram_Power}}

\begin{proof}
First, we show that $W$ is sample-continuous. Since $W$ is centred and has stationary increments, it follows from \eqref{Eq_Condition_Variogram_Power} that, for all $\Mb{x}_1, \Mb{x}_2 \in \mathds{R}^d$, 
\begin{align}
\mathds{E}[(W(\Mb{x}_1)-W(\Mb{x}_2))^2]=\mathrm{Var}(W(\Mb{x}_1)-W(\Mb{x}_2))&=\mathrm{Var}(W(\Mb{x}_1-\Mb{x}_2)-W(\Mb{0})) \nonumber
\\&=\sigma_W^2(\Mb{x}_1-\Mb{x}_2) \nonumber
\\&=\eta \| \Mb{x}_2 - \Mb{x}_1 \|^{\alpha}.
\label{Eq_Expectation_Variogram}
\end{align}
Let us take $C \in (0, \infty)$ and $\rho >0$. As $\alpha >0$, it is well-known that $\lim_{h \to 0} h^{\alpha} |\log(h)|^{1+\rho}=0$. Therefore, there exists $\xi_1>0$ such that, for all $\Mb{x}_1, \Mb{x}_2 \in \mathds{R}^d$ satisfying $\| \Mb{x}_1 - \Mb{x}_2 \|<\xi_1$,
$$ \eta \| \Mb{x}_2 - \Mb{x}_1 \|^{\alpha} | \log(\| \Mb{x}_2 - \Mb{x}_1 \|)|^{1+\rho} \leq C.$$
This means, using \eqref{Eq_Expectation_Variogram}, that there exist $C \in (0, \infty)$ and $\rho, \xi_1 >0$ such that, for all $\Mb{x}_1, \Mb{x}_2 \in \mathds{R}^d$ satisfying $\| \Mb{x}_1 - \Mb{x}_2 \| < \xi_1$,
$$ \mathds{E}[(W(\Mb{x}_1)-W(\Mb{x}_2))^2] \leq \frac{C}{| \log(\| \Mb{x}_2 - \Mb{x}_1 \|)|^{1+\rho}}.$$
Hence, Theorem 1.4.1 in \cite{adler2007random} gives that $W$ is sample-continuous and a.s. bounded on $[0,1]^d$.

\medskip

Second, we prove that the $\sigma_W^2$ satisfies \eqref{Eq_Condition1_Variogram_TCL}.
It follows from \eqref{Eq_Condition_Variogram_Power}, that, for all $\Mb{x} \in [0,1]^d$ and $\Mb{h} \in \mathds{Z}^d$, 
\Beq
\label{Eq_Expression_Difference_Sigma2}
\sigma_W^2(\Mb{h})-\sigma_W^2(\Mb{x}+\Mb{h}) = \eta \left( \| \Mb{h} \|^{\alpha}- \| \Mb{x}+ \Mb{h} \|^{\alpha} \right),
\Eeq
where $\eta>0$ and $\alpha \in (0,2]$. Now, it is well-known that $ \| \Mb{x}+ \Mb{h} \| \geq \| \Mb{h} \| - \| \Mb{x} \|$. This gives, for all $\Mb{h} \in \mathds{Z}^d$ satisfying $\| \Mb{h} \| > d^{\frac{1}{2}}$, that
\Beq
\label{Majoration_Sup_Quantity_Interest}
\sup_{\Mb{x} \in [0,1]^d} \left \{ \| \Mb{h} \|^{\alpha} - \| \Mb{x}+ \Mb{h} \|^{\alpha} \right \} \leq  \sup_{\Mb{x} \in [0,1]^d} \left \{ \| \Mb{h} \|^{\alpha} -   (\| \Mb{h} \| - \| \Mb{x} \|)^{\alpha} \right \}.
\Eeq
Now, we have, for all $\Mb{x} \in [0,1]^d$ and $\Mb{h} \in \mathds{Z}^d$ satisfying $\| \Mb{h} \| > d^{\frac{1}{2}}$, that
$$ (\| \Mb{h} \| - \| \Mb{x} \|)^{\alpha} = \| \Mb{h} \|^{\alpha} \left( 1- \frac{\| \Mb{x} \|}{\| \Mb{h} \|} \right)^{\alpha}.$$
Hence, using a classical Taylor expansion, we obtain that
$$ \| \Mb{h} \|^{\alpha} - (\| \Mb{h} \| - \| \Mb{x} \|)^{\alpha} \underset{\| \Mb{h} \| \to \infty}{=} \alpha \| \Mb{x} \| \| \Mb{h} \|^{\alpha-1} + o \left( \| \Mb{x} \|^2 \| \Mb{h} \|^{\alpha-2} \right).$$ 
Therefore, there exists $B>0$ such that, for all $\Mb{x} \in [0,1]^d$ and $\Mb{h} \in \mathds{Z}^d$ satisfying $\| \Mb{h} \| > d^{\frac{1}{2}}$,
$$ 0 \leq \| \Mb{h} \|^{\alpha} - (\| \Mb{h} \| - \| \Mb{x} \|)^{\alpha} \leq B \| \Mb{x} \| \| \Mb{h} \|^{\alpha-1}.$$
Thus, since $\| \Mb{x} \|$ is bounded on $[0,1]^d$, we directly obtain that
$$
\sup_{\Mb{x} \in [0,1]^d} \left \{ \| \Mb{h} \|^{\alpha} -   (\| \Mb{h} \| - \| \Mb{x} \|)^{\alpha} \right \} \underset{\| \Mb{h} \| \to \infty}{=}o(\| \Mb{h} \|^{\alpha}).$$
Now, we have that 
\Beq
\label{Positivity_Sup_Interest}
\sup_{\Mb{x} \in [0,1]^d} \left \{ \| \Mb{h} \|^{\alpha}- \| \Mb{x}+ \Mb{h} \|^{\alpha} \right \} \geq 0.
\Eeq
Combining \eqref{Eq_Condition_Variogram_Power}, \eqref{Eq_Expression_Difference_Sigma2}, \eqref{Majoration_Sup_Quantity_Interest} and \eqref{Positivity_Sup_Interest}, we obtain
$$\sup_{\Mb{x} \in [0,1]^d} \{ \sigma_W^2(\Mb{h})-\sigma_W^2(\Mb{x}+\Mb{h}) \}\underset{\| \Mb{h} \| \to \infty}{=}o(\sigma_W^2(\Mb{h})).$$
Hence, \eqref{Eq_Condition1_Variogram_TCL} is satisfied.

\medskip

Third, it is clear, using \eqref{Eq_Condition_Variogram_Power}, that 
$$\lim_{\| \Mb{h} \| \to \infty} \frac{\ln(\| \Mb{h} \|)}{\sigma^2_W(\Mb{h})}= 0.$$
Hence, $\sigma_W^2$ satisfies \eqref{Eq_Condition2_Variogram_TCL}.

Finally, since $W$ is a.s. bounded on $[0,1]^d$, Proposition \ref{Prop_Upper_Bound_Expectation_Min_Brown_Resnick} yields that the random field $Y$ satisfies Condition \eqref{Eq_Decay_Cubes_Y} for all $b>0$.
\end{proof}

\subsection{For Theorem \ref{Th_2_TCL_Brown_Resnick}}

\begin{proof}
We assume that $Z$ has been built with a Gaussian random field $\{ W(\Mb{x}) \}_{\Mb{x} \in \mathds{R}^d}$ having variogram $\gamma_W$. 
We consider the random field $\{ W_1(\Mb{x}) \}_{\Mb{x} \in \mathds{R}^d}=\{ W(\Mb{x})-W(\Mb{0}) \}_{\Mb{x} \in \mathds{R}^d}$ and we denote by $Z_1$ the Brown-Resnick random field built with $W_1$. It is clear that $W_1$ is a centred Gaussian random field with stationary increments such that $W_1(\Mb{0})=0$. The random fields $W$ and $W_1$ have the same variogram, and thus the variogram of $W_1$ is written $\gamma_{W_1}(\Mb{x})=\eta \| \Mb{x} \|^{\alpha}$, $\Mb{x} \in \mathds{R}^d$, where $\eta >0$ and $\alpha \in (0,2]$. Now, $\gamma_{W_1}(\Mb{x})=\mathrm{Var}(W_1(\Mb{x})-W_1(\Mb{0}))=\mathrm{Var}(W_1(\Mb{x}))=\sigma_{W_1}^2(\Mb{x})$. Hence, for all $\Mb{x} \in \mathds{R}^d$, $\sigma_{W_1}^2(\Mb{x})=\eta \| \Mb{x} \|^{\alpha}$. Thus, it follows from Proposition \ref{Prop_Variogram_Power} that $W_1$ is sample-continuous, which directly gives that $W$ is sample-continuous. Therefore, applying Proposition 13 in \cite{kabluchko2009stationary}, we obtain that $Z$ and $Z_1$ are sample-continuous. As $W$ and $W_1$ have the same variogram, $Z$ and $Z_1$ have the same finite-dimensional distributions. Moreover, since $Z$ and $Z_1$ are sample-continuous, they have the same distribution in the sense of the induced measure on the space of continuous functions from $\mathds{R}^d$ to $(0, \infty)$. Consequently, we can assume that $Z$ has been built with $W_1$. 

The random field $Z$ is simple (by definition), stationary \citep[see][Theorem 2]{kabluchko2009stationary} and sample-continuous. Moreover, Proposition \ref{Prop_Variogram_Power} gives that the random field $Y$ defined by 
$$ \{ Y(\Mb{x}) \}_{\Mb{x} \in \mathds{R}^d}=\left \{ \exp \left( W_1(\Mb{x})-\frac{\sigma_{W_1}^2(\Mb{x})}{2} \right) \right \}_{\Mb{x} \in \mathds{R}^d}$$ 
satisfies Condition \eqref{Eq_Decay_Cubes_Y} for all $b>0$. Hence, Theorem \ref{Th_General_Case} yields the result.
\end{proof}

\subsection{For Theorem \ref{Th_TCL_Smith}}

\begin{proof}
It is known \cite[see, e.g.,][]{huser2013composite} that the Smith random field with covariance matrix $\Sigma$ corresponds to the Brown-Resnick random field associated with the variogram $\gamma(\Mb{x})=\Mb{x}^{'} \Sigma^{-1} \Mb{x}$, $\Mb{x} \in \mathds{R}^d$. This variogram can be rewritten as $\gamma(\Mb{x})= \| \Mb{x} \|_{\Sigma}^2$, where $\| . \|_{\Sigma}$ is the norm associated with the inner product induced by the matrix $\Sigma^{-1}$. Let $\{ W(\Mb{x}) \}_{\Mb{x} \in \mathds{R}^d}$ be a centred Gaussian random field with stationary increments such that $W(\Mb{0})=0$ and $\sigma_W^2(\Mb{x})=\gamma(\Mb{x})$, $\Mb{x} \in \mathds{R}^d$. Since all norms are equivalent in $\mathds{R}^d$, there exist $C_3, C_4>0$ such that, for all $\Mb{x} \in \mathds{R}^d$, $C_3 \| \Mb{x} \| \leq \| \Mb{x} \|_{\Sigma} \leq C_4 \| \Mb{x} \|$. Hence, a hardly modified version of the proof of Proposition \ref{Prop_Variogram_Power} leads that $W$ is sample-continuous and that the random field
$$ \{ Y(\Mb{x}) \}_{\Mb{x} \in \mathds{R}^d}=\left \{ \exp \left( W(\Mb{x})-\frac{\sigma_W^2(\Mb{x})}{2} \right) \right \}_{\Mb{x} \in \mathds{R}^d}$$
satisfies Condition \eqref{Eq_Decay_Cubes_Y} for all $b>0$. Then, the same proof as that of Theorem \ref{Th_2_TCL_Brown_Resnick} leads to the result.
\end{proof}

\newpage
\bibliographystyle{apalike}
\bibliography{../../../../../Bibliography_Erwan/References_Erwan}

\end{document}